\documentclass[11pt]{article}

\usepackage[a4paper,margin=1in]{geometry}
\usepackage{amsmath,amssymb,amsthm,mathtools,bm}
\usepackage{mathrsfs}
\usepackage{hyperref}
\hypersetup{hidelinks}
\usepackage{xcolor}
\usepackage{graphicx}
\usepackage{epstopdf}
\usepackage{placeins}
\usepackage{booktabs}
\usepackage{multirow}

\theoremstyle{plain}

\numberwithin{equation}{section}

\newtheorem{proposition}{Proposition}[section]

\newtheorem{corollary}[proposition]{Corollary}
\newtheorem{remark}{Remark}[section]

\renewcommand{\arraystretch}{1.25}

\title{A decoupled structure-preserving scheme for pressure-driven multispecies tumor growth}

\author{
Weijie Huang\thanks{School of Mathematics and Statistics, Beijing Jiaotong
University, Beijing 100044, People's Republic of China; Beijing Key
Laboratory of Biological Big Data and Topological Statistics, Beijing
Jiaotong University, Beijing 100044, People's Republic of China
{\tt (wjhuang@bjtu.edu.cn)}.}
\and
Xinran Ruan\thanks{(Corresponding author) School of Mathematical Sciences,
Capital Normal University, Beijing 100048, People's Republic of China
{\tt (xinran.ruan@cnu.edu.cn)}.}
}
\date{}

\begin{document}

\maketitle

\begin{abstract}
We propose an Onsager-based scheme for an $m$-species
pressure-driven reaction--cross-diffusion system with a common
homeostatic pressure.
The scheme uses a modified scalar energy structure to update the total density, followed by compatible
updates of the individual species.
The global nonlinear problem remains scalar, and the species densities
sum exactly to the updated total density.
We establish unique solvability and show that the scheme preserves
species nonnegativity, the pressure bound and total-mass monotonicity
and satisfies a discrete energy dissipation law without a time-step
restriction.
We further develop a predictor--corrector extension to improve temporal
accuracy while retaining these structural properties.
Numerical experiments verify the accuracy and structure-preserving
properties of the proposed schemes and demonstrate their computational
efficiency.
\end{abstract}



\section{Introduction}

Mechanical interactions are increasingly recognized as important regulators
of solid-tumor growth and progression
\cite{ButcherAllistonWeaver2009,NiaMunnJain2020}.
As tumor cells proliferate in a confined environment, competition for space
and resistance from the surrounding tissue generate mechanical stresses
that can in turn regulate tissue growth and cell proliferation
\cite{Shraiman2005}.
Experiments on multicellular tumor spheroids have directly shown that
compression can substantially suppress tumor growth and proliferation
\cite{DelarueEtAl2014, MontelEtAl2011}.
The concept of homeostatic pressure, at which cell proliferation and loss
balance, provides a useful mechanical description of tissue growth and
competition
\cite{BasanRislerJoannySastreGarauProst2009}.
At the same time, tumors are highly heterogeneous, containing cell
populations with distinct phenotypic and functional properties that arise
from genetic differences, environmental variation, and reversible changes
in cellular states
\cite{MarusykAlmendroPolyak2012,MeachamMorrison2013}.
For spatial tumor models, these mechanical and biological effects make it
important to resolve not only overall tumor growth but also the evolution
of its internal cellular composition.

Pressure-driven continuum models provide a mathematical framework for
describing the mechanical regulation of tissue growth. For a single cell
population, density-dependent pressure together with Darcy transport and
pressure-regulated proliferation leads to degenerate
porous-medium-type equations
\cite{PerthameQuirosVazquez2014,PerthameVauchelet2015}.
A substantial mathematical theory has been developed for these models.
In particular, the stiff-pressure limit connects the compressible density
formulation to incompressible Hele--Shaw-type free-boundary problems
\cite{DavidPerthame2021,KimPozar2018}.

The pressure-driven framework has subsequently been extended to
heterogeneous cell populations. A two-species hyperbolic--parabolic
model was analyzed in
\cite{GwiazdaPerthameSwierczewskaGwiazda2019},
and incompressible limits for related two-species systems were established
in
\cite{BubbaPerthamePoucholSchmidtchen2020,DegondHechtVauchelet2020}.
Phenotype-structured and reaction--cross-diffusion models have also been
considered in
\cite{David2023PhenotypicHeterogeneity,DouLiuZhou2023},
while more detailed multispecies models have been used to describe
compositional heterogeneity in tumor evolution
\cite{FrieboesJinChuangWiseLowengrubCristini2010,WiseLowengrubFrieboesCristini2008}.

In this paper, we consider an $m$-species pressure-driven
reaction--cross-diffusion system that generalizes the
two-population model in \cite{BubbaPerthamePoucholSchmidtchen2020},
\begin{equation}
\begin{cases}
\displaystyle
\partial_t \boldsymbol n
-\nabla\cdot(\boldsymbol n\nabla p)
=\mathbf R(p)^T\boldsymbol n,\\[1ex]
\displaystyle
p=n^\gamma,
\qquad
n=\mathbf 1^T\boldsymbol n,
\qquad \gamma>1.
\end{cases}
\label{eq:intro-model}
\end{equation}
Here $\boldsymbol n=(n_1,\ldots,n_m)^T$ is the
species-density vector, where $n_i=n_i(x,t)\geq0$
denotes the density of species $i$,  $\mathbf 1=(1,\ldots,1)^T\in\mathbb R^m$,
$n$ is the total density, and $p$ is the mechanical
pressure generated by cell crowding.
All species are transported by the same Darcy velocity
$-\nabla p$.
The reaction matrix is given by
$\mathbf R(p)=(R_{\ell i}(p))_{\ell,i=1}^m$,
where $R_{\ell i}\in C^1([0,\infty))$ describes
the contribution of species $\ell$ to species $i$
through growth, loss, and transitions.
On a bounded domain \(\Omega\subset\mathbb R^d\), we impose
\begin{equation}\label{eq:boundary_condition}
  n_i\nabla p\cdot\nu=0
  \quad\text{on }\partial\Omega,
  \qquad
  n_i(\cdot,0)=n_{i,0},
  \qquad
  i=1,\ldots,m,
\end{equation}
where \(\nu\) is the outward unit normal to \(\partial\Omega\).

Multiplying the first equation in \eqref{eq:intro-model}  by \(\mathbf1^T\) gives
\begin{equation}
 \label{eq:intro-total-equation}
  \partial_t n-\nabla\cdot(n\nabla p)
  =
  \sum_{\ell=1}^m n_\ell H_\ell(p),
  \quad \text{ where } \quad
  H_\ell(p):=\sum_{i=1}^m R_{\ell i}(p).
\end{equation}
Here \(H_\ell\) is the net growth rate associated with species \(\ell\).
We assume that all species share a common
homeostatic pressure \(p_H>0\),
\begin{equation}
  H_\ell(p_H)=0,
  \qquad
  H_\ell'(p)\leq0,
  \qquad
  p\geq0,
  \qquad
  \ell=1,\ldots,m.
  \label{eq:intro-common-ph}
\end{equation}
Hence each species makes a nonnegative contribution to total
growth below $p_H$ and a nonpositive contribution above $p_H$,
while its growth and transition rates may differ from those
of the other species.
The common homeostatic-pressure assumption is consistent
with pressure-inhibited growth
\cite{BasanRislerJoannySastreGarauProst2009,DelarueEtAl2014,MontelEtAl2011,Shraiman2005}.
We further assume in
\eqref{eq:intro-model} that
  $R_{\ell i}(p)\geq0$,
 for
  $\ell\ne i \text{ and } p\in[0,p_H]$.

The numerical approximation of pressure-driven growth and related
degenerate diffusion equations has been studied from several perspectives.
Bessemoulin-Chatard and Filbet~\cite{BessemoulinChatardFilbet2012}
proposed a finite-volume scheme for nonlinear degenerate parabolic
equations, while Duan et al.~\cite{DuanLiuWangYue2019} developed
numerical methods for the porous-medium equation based on an energetic
variational approach. For a related reaction--cross-diffusion system,
Carrillo et al.~\cite{CarrilloFagioliSantambrogioSchmidtchen2018}
introduced a variational splitting scheme and studied the preservation of
segregation in the presence of vacuum. For the single-species
pressure-driven tumor model, Liu et al.~\cite{LiuTangWangZhou2018}
developed a front-capturing scheme and later considered models coupled to
nutrients~\cite{LiuTangWangZhou2019}. David and
Ruan~\cite{DavidRuan2022} proposed an asymptotic-preserving scheme, while
Huang and Ruan~\cite{HuangRuan2026} derived an Onsager variational
scheme.
These studies provide important numerical tools for related degenerate and
pressure-driven models. However, numerical methods that simultaneously
exploit the structure of the total-density dynamics and resolve the
evolution of individual species remain less developed.

For the multispecies system \eqref{eq:intro-model}, the numerical
challenge is to preserve the structure of the total-density dynamics
while resolving the individual species consistently.
To address this challenge, we identify a scalar energy structure for the
total-density dynamics and build on the Onsager formulation of
\cite{HuangRuan2026} to construct a first-order multispecies scheme.
The scheme couples a scalar total-density update with species updates
that reconstruct the updated total density exactly.
We establish unique solvability and show that the scheme preserves
nonnegativity, the pressure bound and total-mass monotonicity and satisfies
a discrete energy dissipation law without a time-step restriction.
We further develop a predictor--corrector extension to improve temporal
accuracy while retaining these structural properties.
Numerical experiments demonstrate the accuracy and computational
performance of both schemes.

The remainder of the paper is organized as follows.
Section~\ref{sec:multispecies-model} introduces the composition variables,
discusses the non-closure of the total-density equation, and derives the
scalar energy structure. Section~\ref{sec:multi-scheme} presents the
numerical scheme and establishes its structure-preserving
properties. Section~\ref{sec:numerical} contains the numerical
experiments. The two-dimensional extension is given in
Appendix~\ref{app:2d}.

\section{Model structure}
\label{sec:multispecies-model}

In this section, we present the structural properties of the multispecies
model that will be used to construct the numerical scheme.

\subsection{Total density and composition}
\label{subsec:multi-composition}
\label{subsec:nonclosure}

To separate the evolution of the total density from that
of the species composition, we introduce the local
fractions of each species. On the occupied set
$\{n>0\}$, define
$\phi_i=\frac{n_i}{n}$,
$i=1,\ldots,m$.
The composition vector
$\boldsymbol\phi=(\phi_1,\ldots,\phi_m)^T$ belongs to
\begin{equation}\label{eq:simplex}
  \Delta_m
  =
  \left\{
    \boldsymbol\phi\in\mathbb R^m
    \mathrel{}\middle|\mathrel{}
    \phi_i\geq0,\quad
    \sum_{i=1}^m\phi_i=1
  \right\}.
\end{equation}
Using the net growth rates $H_\ell$ defined in
\eqref{eq:intro-total-equation}, let
  $G(\boldsymbol\phi,p)
  =
  \sum_{\ell=1}^m\phi_\ell H_\ell(p)$.
The total-density equation \eqref{eq:intro-total-equation}
then becomes
\begin{equation}\label{eq:multi-density-phi}
  \partial_t n-\nabla\cdot(n\nabla p)
  =
  nG(\boldsymbol\phi,p),
  \qquad p=n^\gamma.
\end{equation}
It is remarkable that, unless the net growth rates $H_\ell$ coincide, the total
growth rate depends on the local composition, and
\eqref{eq:multi-density-phi} is not closed in $n$ and $p$.

Substituting $n_i=n\phi_i$ into the species equations
\eqref{eq:intro-model} and using
\eqref{eq:multi-density-phi}, we obtain
\begin{equation}\label{eq:composition-equation}
  \partial_t\phi_i-\nabla p\cdot\nabla\phi_i
  =
  B_i(\boldsymbol\phi,p),
  \qquad i=1,\ldots,m,
\end{equation}
where
\begin{equation}\label{eq:composition-reaction}
  B_i(\boldsymbol\phi,p)
  =
  \sum_{\ell=1}^m\phi_\ell R_{\ell i}(p)
  -
  \phi_iG(\boldsymbol\phi,p).
\end{equation}
By the definitions of $H_\ell$ and $G$,
$\sum_{i=1}^m B_i(\boldsymbol\phi,p)=0$,
consistent with $\sum_{i=1}^m\phi_i=1$.

The composition variables are undefined on $\{n=0\}$,
where all species densities vanish and the total reaction
term is zero.

\subsection{Basic properties}
\label{subsec:invariant-region}

We establish the nonnegativity of the species densities,
the pressure bound and total-mass monotonicity.

\begin{proposition}[Basic bounds]
\label{prop:basic-bounds}
Suppose that the initial densities satisfy
$n_{i,0}\geq0$ and $0\leq p_0\leq p_H$, where
$
  p_0=\left(\sum_{i=1}^m n_{i,0}\right)^\gamma.
$
Let $(\boldsymbol n,p)$ be a sufficiently smooth solution
of \eqref{eq:intro-model} satisfying the no-flux boundary
condition \eqref{eq:boundary_condition}. Then, as long as
the solution exists,
\begin{equation}\label{eq:continuous-basic-bounds}
  n_i\geq0,\qquad
  0\leq p\leq p_H,\qquad
  0\leq n\leq n_H:=p_H^{1/\gamma}.
\end{equation}
Moreover, $\boldsymbol\phi\in\Delta_m$ on $\{n>0\}$,
and the total mass satisfies
\begin{equation}\label{eq:continuous-mass-estimate}
  \frac{d}{dt}\int_\Omega n\,dx\geq0,
  \qquad
  \int_\Omega n(x,0)\,dx
  \leq\int_\Omega n(x,t)\,dx
  \leq n_H|\Omega|.
\end{equation}
\end{proposition}

\begin{proof}
Multiplying \eqref{eq:intro-total-equation} by
$\gamma n^{\gamma-1}$ and using $p=n^\gamma$, we obtain
\begin{equation}\label{eq:pressure-equation-for-bound}
  \partial_t p
  =
  |\nabla p|^2+\gamma p\Delta p
  +\gamma n^{\gamma-1}
   \sum_{\ell=1}^m n_\ell H_\ell(p).
\end{equation}
Since $H_\ell(p_H)=0$, the pressure bound follows
from the standard comparison argument.

To prove species nonnegativity, we consider the characteristics
$\dot{\boldsymbol\xi}=-\nabla p$, along which
\[
  \frac{d}{dt}\boldsymbol n(\boldsymbol\xi(t),t)
  =
  \left(\Delta p\,\mathbf I+\mathbf R(p)^T\right)
  \boldsymbol n(\boldsymbol\xi(t),t).
\]
The nonnegative off-diagonal entries of $\mathbf R(p)$
ensure that nonnegative initial densities remain
nonnegative. Consequently, $n\geq0$ and
$\boldsymbol\phi\in\Delta_m$ wherever $n>0$.

Finally, \eqref{eq:intro-common-ph} gives
$H_\ell(p)\geq0$ for $0\leq p\leq p_H$.
Since $n_\ell\geq0$, integrating
\eqref{eq:intro-total-equation} over $\Omega$
and using the no-flux boundary condition yields
total-mass monotonicity.
The upper bound follows from $n\leq n_H$.
\end{proof}

\subsection{Modified energy and Onsager structure}
\label{subsec:multi-modified-energy}

We next examine the energy dissipation and Onsager structure
of the total-density equation.
Following \cite{HuangRuan2026}, we define the modified energy
and its variational derivative by
\begin{equation}\label{eq:mu}
  E(n)=\int_\Omega
  \left(\frac{n^{\gamma+1}}{\gamma+1}-p_H n\right)\,dx,
  \qquad
  \mu=\frac{\delta E}{\delta n}=p-p_H.
\end{equation}
Let
\begin{equation}\label{eq:continuous-secant-coefficient}
  \sigma_\ell(p)=
  \begin{cases}
    \dfrac{H_\ell(p)}{p-p_H}, & p\ne p_H,\\[1ex]
    H_\ell'(p_H), & p=p_H.
  \end{cases}
\end{equation}
By \eqref{eq:intro-common-ph}, $\sigma_\ell$ is continuous
and nonpositive, so the reaction mobility and total reaction
term satisfy
\begin{align}
  M(\boldsymbol n,p)
  &=-\sum_{\ell=1}^m n_\ell\sigma_\ell(p)\geq0, \quad
  \sum_{\ell=1}^m n_\ell H_\ell(p)
  =-M(\boldsymbol n,p)\mu.
  \label{eq:multi-reaction-mobility}
\end{align}
For smooth nonnegative solutions, testing
\eqref{eq:intro-total-equation} with $\mu$
and using the no-flux boundary condition gives
\begin{equation}\label{eq:multi-onsager-energy-law}
  \frac{d}{dt}E(n)
  =\int_\Omega\mu\partial_t n\,dx
  =-\int_\Omega n|\nabla p|^2\,dx
   -\int_\Omega M|p-p_H|^2\,dx
  \leq0.
\end{equation}

The total-density equation \eqref{eq:intro-total-equation}
can also be expressed in an Onsager variational form
\cite{ChenLiuXu2025,Doi2013,Onsager1931}.
At a given state $(\boldsymbol n,p)$, we introduce
the Rayleighian $\mathcal R(\boldsymbol j,r)$,
which combines the energy rate $\frac{d}{dt}E(n)$
with a nonnegative dissipation potential
$\Phi(\boldsymbol j,r)$.
The transport flux $\boldsymbol j$ and reaction
rate $r$ are treated as process variables subject
to the balance law
\begin{equation}\label{eq:balance_law_and_boundary}
  \partial_t n+\nabla\cdot\boldsymbol j=r
  \quad\text{in }\Omega,
  \qquad
  \boldsymbol j\cdot\boldsymbol\nu=0
  \quad\text{on }\partial\Omega.
\end{equation}
The evolution is then determined by minimizing this
functional $\mathcal R(\boldsymbol j,r)$ with respect
to the process variables $\boldsymbol j,r$ at the
given state $(\boldsymbol n,p)$.

As in \cite{HuangRuan2026}, we consider the quadratic
dissipation potential
\begin{equation}\label{eq:dissipation}
  \Phi(\boldsymbol j,r)
  =
  \frac12\int_\Omega\frac{|\boldsymbol j|^2}{n}\,dx
  +
  \frac12\int_\Omega\frac{r^2}{M}\,dx.
\end{equation}
Here $M$ is the reaction mobility in \eqref{eq:multi-reaction-mobility}.
The quotients are understood in
the extended-value sense, requiring
$\boldsymbol j=\boldsymbol0$ on $\{n=0\}$ and
$r=0$ on $\{M=0\}$.
The corresponding Rayleighian is
\begin{equation}\label{eq:multi-rayleighian}
\begin{aligned}
  \mathcal R(\boldsymbol j,r)
  &=
  \Phi(\boldsymbol j,r)+\frac{d}{dt}E(n)
  =
  \Phi(\boldsymbol j,r)
  +\int_\Omega\nabla\mu\cdot\boldsymbol j\,dx
  +\int_\Omega\mu r\,dx.
\end{aligned}
\end{equation}
Minimizing $\mathcal R$ with respect to
$(\boldsymbol j,r)$ at the fixed state yields
\[
  \boldsymbol j=-n\nabla\mu=-n\nabla p,
  \qquad
  r=-M\mu
  =\sum_{\ell=1}^m n_\ell H_\ell(p).
\]
Substituting these relations into the balance law
recovers the total-density equation
\eqref{eq:intro-total-equation}.
The Onsager formulation provides the basis for the
total-density update in our numerical scheme.

\section{Numerical scheme and structure preservation}
\label{sec:multi-scheme}

We construct a numerical scheme that combines a scalar
density update with the evolution of individual species.
The scalar step determines the total reaction rate
and transport flux, which are then used to update the
species densities consistently.

\subsection{Scalar step for the total density}
\label{subsec:multi-time-discrete}

We first discretize the total-density equation using
the Onsager formulation established in
Section~\ref{subsec:multi-modified-energy}.
Let $\tau>0$ be the time step and $t^k=k\tau$.
Given the species densities
$\boldsymbol n^k=(n_1^k,\ldots,n_m^k)^T$,
we set
\[
  n^k=\mathbf1^T\boldsymbol n^k,
  \qquad
  p^k=(n^k)^\gamma,
  \qquad
  M^k
  =
  -\sum_{\ell=1}^m n_\ell^k\sigma_\ell(p^k)
  \geq0,
\]
where $\mathbf1=(1,\ldots,1)^T\in\mathbb R^m$
and the mobility $M^k$ follows \eqref{eq:multi-reaction-mobility}.
%

Following the time discretization in
\cite{HuangRuan2026}, we replace the energy rate
by the energy increment over one time step and
evaluate the transport and reaction mobilities
at time $t^k$.
This gives the discrete Rayleighian
\begin{equation}\label{eq:multi-discrete-rayleighian}
  \mathcal R_\tau^k(n,\boldsymbol j,r)
  =
  E(n)-E(n^k)
  +\frac{\tau}{2}\int_\Omega
    \frac{|\boldsymbol j|^2}{n^k}\,dx
  +\frac{\tau}{2}\int_\Omega
    \frac{r^2}{M^k}\,dx.
\end{equation}
As in the continuous formulation, the quotients
are understood in the extended-value sense
when the mobilities vanish.
The updated density and the process variables are
determined by minimizing $\mathcal R_\tau^k$ subject to
\begin{equation}\label{eq:multi-discrete-balance}
  n-n^k+\tau\nabla\cdot\boldsymbol j=\tau r,
  \qquad
  n\geq0,\qquad
  \boldsymbol j\cdot\boldsymbol\nu=0
  \quad\text{on }\partial\Omega.
\end{equation}
The optimality conditions yield
$
  \boldsymbol j^{k+1}
  =-n^k\nabla p^{k+1}$, $
  r^{k+1}
  =-M^k(p^{k+1}-p_H).
$
Substituting these into
\eqref{eq:multi-discrete-balance}, we obtain
\begin{equation}\label{eq:multi-time-scheme}
\left\{
\begin{aligned}
  &\frac{n^{k+1}-n^k}{\tau}
  =
  \nabla\cdot\bigl(n^k\nabla p^{k+1}\bigr)
  -M^k\bigl(p^{k+1}-p_H\bigr)
  &&\text{in }\Omega,\\
  &p^{k+1}=(n^{k+1})^\gamma
  &&\text{in }\Omega,\\
  &n^k\nabla p^{k+1}\cdot\nu=0
  &&\text{on }\partial\Omega.
\end{aligned}
\right.
\end{equation}

Since $E$ is convex, multiplying the first equation
in \eqref{eq:multi-time-scheme} by
$\mu^{k+1}=p^{k+1}-p_H$ and integrating by parts gives
\begin{equation}\label{eq:multi-semidiscrete-energy}
  E(n^{k+1})
  +\tau\int_\Omega n^k|\nabla p^{k+1}|^2\,dx
  +\tau\int_\Omega M^k|p^{k+1}-p_H|^2\,dx
  \leq E(n^k).
\end{equation}
Thus, the time-discrete scalar step preserves the
modified energy dissipation structure
\eqref{eq:multi-onsager-energy-law}.

\subsection{Fully discrete scheme}
\label{subsec:multi-full-discrete}

We discretize \eqref{eq:multi-time-scheme} using
finite differences in conservative form.
We present the one-dimensional scheme on
$\Omega=(a,b)$, while the two-dimensional extension
is given in Appendix~\ref{app:2d}.
We divide $\Omega$ into $N$ uniform cells of width
$h=(b-a)/N$, with cell centers
$x_j=a+(j-1/2)h$, $j=1,\ldots,N$.
For values $z_j$ at the cell centers and
$q_{j+1/2}$ at the cell interfaces, define
\begin{equation}\label{eq:multi-difference-operators}
  (\delta_h z)_{j+1/2}
  =
  \frac{z_{j+1}-z_j}{h},
  \qquad
  (d_hq)_j
  =
  \frac{q_{j+1/2}-q_{j-1/2}}{h}.
\end{equation}
Unless otherwise stated, cell equations hold for
$j=1,\ldots,N$, and interior interface equations
for $j=1,\ldots,N-1$.

\paragraph{Step 1: Scalar density step.}

We first compute the updated total density and pressure.
Given the species densities
$\boldsymbol n_j^k=(n_{1,j}^k,\ldots,n_{m,j}^k)^T$,
we set
\[
  n_j^k=\mathbf1^T\boldsymbol n_j^k,
  \qquad
  p_j^k=(n_j^k)^\gamma,
  \qquad
  M_j^k
  =-\sum_{\ell=1}^m
  n_{\ell,j}^k\sigma_\ell(p_j^k).
\]
Following \eqref{eq:multi-time-scheme}, we solve
\begin{equation}\label{eq:multi-discrete-density-step}
  \frac{n_j^{k+1}-n_j^k}{\tau}
  =
  (d_hF^{k+1})_j
  -M_j^k(p_j^{k+1}-p_H),
  \qquad
  p_j^{k+1}=(n_j^{k+1})^\gamma,
\end{equation}
where
\begin{equation}\label{eq:multi-discrete-flux}
\begin{gathered}
  F_{j+1/2}^{k+1}
  =n_{j+1/2}^k(\delta_hp^{k+1})_{j+1/2}, \quad
    n_{j+1/2}^k=\frac{n_j^k+n_{j+1}^k}{2},
  \quad
  F_{1/2}^{k+1}=F_{N+1/2}^{k+1}=0.
\end{gathered}
\end{equation}
This gives a nonlinear system for the total density
alone, which may be solved by a damped Newton method.

For the subsequent species updates, we define
\begin{equation}\label{eq:multi-source-local-mass}
  r_j^{k+1}
  =
  -M_j^k(p_j^{k+1}-p_H),
  \qquad
  \widehat n_j^{k+1}
  =
  n_j^k+\tau r_j^{k+1}.
\end{equation}
Here $r_j^{k+1}$ is the total reaction rate, and
$\widehat n_j^{k+1}$ is the total density to be
distributed among the species in Step~2.
The scalar balance can then be written as
\begin{equation}\label{eq:multi-density-conservative-step}
    n_j^{k+1}-\tau(d_hF^{k+1})_j
  =\widehat n_j^{k+1}.
\end{equation}

\paragraph{Step 2: Species reaction step.}

To preserve species nonnegativity and consistency
with the total-density update, we split the species
update into reaction and transport steps.
We first update the local composition and distribute
the total density $\widehat n_j^{k+1}$ computed
in Step~1 among the species.

For $n_j^k>0$, let
$\boldsymbol\phi_j^k=\boldsymbol n_j^k/n_j^k$.
With the pressure fixed at $p_j^{k+1}$,
the reaction part of \eqref{eq:composition-equation}
reduces to the local ODE system
$\dot\phi_i=B_i(\boldsymbol\phi,p_j^{k+1})$,
with $\boldsymbol\phi(0)=\boldsymbol\phi_j^k$.
Here $B_i$ is defined in
\eqref{eq:composition-reaction}.

The solution of the above ODE system is obtained
by normalizing the solution of
$\dot{\boldsymbol y}
=\mathbf R(p_j^{k+1})^T\boldsymbol y$,
with $\boldsymbol y(0)=\boldsymbol\phi_j^k$.
Evaluating at time $\tau$ gives
\begin{equation}\label{eq:multi-reaction-composition}
  \widehat{\boldsymbol\phi}_j^{k+1}
  =
  \frac{
    \exp\!\left(\tau\mathbf R(p_j^{k+1})^T\right)
    \boldsymbol\phi_j^k
  }{
    \mathbf1^T
    \exp\!\left(\tau\mathbf R(p_j^{k+1})^T\right)
    \boldsymbol\phi_j^k
  }.
\end{equation}

If $n_j^k=0$, then
$\boldsymbol n_j^k=\boldsymbol0$ and $M_j^k=0$.
Equation~\eqref{eq:multi-source-local-mass}
therefore gives
$r_j^{k+1}=\widehat n_j^{k+1}=0$,
so no composition is needed at these cells.

We assign the prescribed total density
to the species by
\begin{equation}\label{eq:multi-post-reaction_2}
  \widehat{\boldsymbol n}_j^{k+1}
  =
  \begin{cases}
    \widehat n_j^{k+1}
    \widehat{\boldsymbol\phi}_j^{k+1},
      & n_j^k>0,\\[1ex]
    \boldsymbol0,
      & n_j^k=0.
  \end{cases}
\end{equation}

\paragraph{Step 3: Species transport step.}

We transport the species using the total flux
computed in Step~1, requiring the species fluxes
to sum to the prescribed total flux.
Let $\boldsymbol\phi_j^{k+1}$ denote the unknown
composition after transport, and set
$s^+=\max\{s,0\}$ and $s^-=\min\{s,0\}$.
Since $F^{k+1}$ has the opposite sign to the
physical flux, we use the upwind species flux
\begin{equation}\label{eq:multi-upwind-species-flux}
  \boldsymbol F_{j+1/2}^{k+1}
  =
  (F_{j+1/2}^{k+1})^-
  \boldsymbol\phi_j^{k+1}
  +
  (F_{j+1/2}^{k+1})^+
  \boldsymbol\phi_{j+1}^{k+1},
\end{equation}
with
$\boldsymbol F_{1/2}^{k+1}
=\boldsymbol F_{N+1/2}^{k+1}
=\boldsymbol0$.
The compositions satisfy
\begin{equation}\label{eq:multi-species-conservative-step}
  n_j^{k+1}\boldsymbol\phi_j^{k+1}
  -\tau(d_h\boldsymbol F^{k+1})_j
  =\widehat{\boldsymbol n}_j^{k+1}.
\end{equation}

To solve this system, define
\begin{equation}\label{eq:multi-active-diagonal}
  D_j^{k+1}
  =
  n_j^{k+1}
  +
  \frac{\tau}{h}
  \left[
    (F_{j-1/2}^{k+1})^+
    -
    (F_{j+1/2}^{k+1})^-
  \right],
\end{equation}
and the active set
$
  \mathcal A^{k+1}
  =
  \{j:D_j^{k+1}>0\}.
$
For $j\in\mathcal A^{k+1}$,
\eqref{eq:multi-species-conservative-step} becomes
\begin{equation}\label{eq:multi-remapping-recursion}
  D_j^{k+1}\boldsymbol\phi_j^{k+1}
  =
  \widehat{\boldsymbol n}_j^{k+1}
  +\frac{\tau}{h}
  \left[
    (F_{j+1/2}^{k+1})^+
      \boldsymbol\phi_{j+1}^{k+1}
    -
    (F_{j-1/2}^{k+1})^-
      \boldsymbol\phi_{j-1}^{k+1}
  \right].
\end{equation}
By \eqref{eq:multi-discrete-flux},
a neighboring composition has a nonzero coefficient
only when that cell has higher pressure,
in which case it also belongs to $\mathcal A^{k+1}$.
We therefore sort the active cells by
$p_j^{k+1}$ from largest to smallest and solve
\eqref{eq:multi-remapping-recursion} successively.
At each cell, the right-hand side is already known,
and division by $D_j^{k+1}>0$ gives
$\boldsymbol\phi_j^{k+1}$.
Cells with equal pressure have zero flux between them
and may be processed in any order.

For $j\notin\mathcal A^{k+1}$,
\eqref{eq:multi-active-diagonal} and the scalar balance
\eqref{eq:multi-density-conservative-step},
together with nonnegativity, imply
\[
  n_j^{k+1}=\widehat n_j^{k+1}=0,
  \qquad
  F_{j-1/2}^{k+1}=F_{j+1/2}^{k+1}=0.
\]
The post-reaction species densities also vanish,
so no composition is needed at these cells.

Finally, we reconstruct the species densities by
\begin{equation}\label{eq:multi-species-reconstruction}
  \boldsymbol n_j^{k+1}
  =
  \begin{cases}
    n_j^{k+1}\boldsymbol\phi_j^{k+1},
      & j\in\mathcal A^{k+1},\\
    \boldsymbol0,
      & j\notin\mathcal A^{k+1}.
  \end{cases}
\end{equation}

\subsection{Structure-preserving properties}
\label{subsec:multi-discrete-structure}

We establish the nonnegativity, pressure bound,
and energy dissipation of the fully discrete scheme,
together with the compatibility between the species
and scalar updates.

Define the discrete energy by
\begin{equation}\label{eq:multi-discrete-energy}
  E_h(n)=h\sum_{j=1}^N\left(\frac{1}{\gamma+1}n_j^{\gamma+1}-p_Hn_j\right).
\end{equation}

\begin{proposition}[Scalar step]
\label{prop:multi-scalar-step}
Suppose that, at time level $k$,
$
  n_{i,j}^k\geq0$, $n_j^k=\sum_{i=1}^m n_{i,j}^k$, $0\leq p_j^k=(n_j^k)^\gamma\leq p_H.
$
Then \eqref{eq:multi-discrete-density-step}
admits a unique solution satisfying
\[
  0\leq n_j^{k+1}\leq n_H,\qquad 0\leq p_j^{k+1}\leq p_H,\qquad r_j^{k+1}\geq0,
\]
where $n_H=p_H^{1/\gamma}$.
Moreover, the total mass satisfies
\[
  h\sum_{j=1}^N n_j^k\leq h\sum_{j=1}^N n_j^{k+1}\leq n_H|\Omega|,
\]
and the discrete energy satisfies
\begin{equation}\label{eq:multi-discrete-energy-ineq}
  E_h(n^{k+1})+\tau h\sum_{j=1}^{N-1}n_{j+1/2}^k|(\delta_hp^{k+1})_{j+1/2}|^2+\tau h\sum_{j=1}^N M_j^k|p_j^{k+1}-p_H|^2\leq E_h(n^k).
\end{equation}
\end{proposition}

\begin{proof}
By \eqref{eq:multi-reaction-mobility},
\(n_{\ell,j}^k\geq0\) and \(\sigma_\ell(p_j^k)\leq0\)
imply \(M_j^k\geq0\).
Following \cite{HuangRuan2026}, consider the functional
on \([0,\infty)^N\)
\[
\begin{aligned}
  \mathcal J_h^k(p)
  ={}&h\sum_{j=1}^N
  \left[
    \frac{\gamma}{\gamma+1}p_j^{(\gamma+1)/\gamma}
    -(n_j^k+\tau M_j^kp_H)p_j
    +\frac{\tau}{2}M_j^kp_j^2
  \right]
  +\frac{\tau h}{2}\sum_{j=1}^{N-1}
    n_{j+1/2}^k|(\delta_hp)_{j+1/2}|^2.
\end{aligned}
\]
Coercivity and strict convexity give a unique minimizer
\(p^{k+1}\geq0\).
Set \(n_j^{k+1}=(p_j^{k+1})^{1/\gamma}\).
At cells where \(p_j^{k+1}>0\), the optimality condition
gives \eqref{eq:multi-discrete-density-step}.
If \(p_j^{k+1}=0\), only nonnegative perturbations
of this component are admissible.
Moreover, the neighboring pressures and transport
mobilities are nonnegative, so
\((d_hF^{k+1})_j\geq0\).
The one-sided optimality condition therefore yields
\[
  0\leq
  \frac{\partial\mathcal J_h^k}{\partial p_j}(p^{k+1})
  =
  -h\left[
    n_j^k+\tau M_j^kp_H+\tau(d_hF^{k+1})_j
  \right]
  \leq0.
\]
Thus the derivative vanishes at these cells as well,
and the minimizer solves the discrete system everywhere.
Conversely, every nonnegative solution satisfies
the optimality conditions of \(\mathcal J_h^k\),
so strict convexity gives uniqueness.

To prove the pressure bound, choose a cell \(j_*\)
where \(p^{k+1}\) attains its maximum.
Since all neighboring pressures are no larger,
the nonnegative transport mobilities and zero boundary
fluxes imply \((d_hF^{k+1})_{j_*}\leq0\).
If \(p_{j_*}^{k+1}>p_H\), the reaction term
\(M_{j_*}^k(p_H-p_{j_*}^{k+1})\) is also nonpositive.
Equation \eqref{eq:multi-discrete-density-step}
then gives \(n_{j_*}^{k+1}\leq n_{j_*}^k\leq n_H\),
contradicting
\(n_{j_*}^{k+1}=(p_{j_*}^{k+1})^{1/\gamma}>n_H\).
Hence \(p_j^{k+1}\leq p_H\) in every cell,
and the density bound and \(r_j^{k+1}\geq0\) follow.
Summing \eqref{eq:multi-discrete-density-step}
and using the zero boundary fluxes shows that
the total mass is nondecreasing.
Together with \(n_j^{k+1}\leq n_H\),
this proves the mass estimate.

Finally, set \(\mu_j^{k+1}=p_j^{k+1}-p_H\).
Convexity of \(E_h\), followed by testing
\eqref{eq:multi-discrete-density-step}
with \(\tau h\mu_j^{k+1}\) and discrete summation
by parts, gives
\[
\begin{aligned}
  E_h(n^{k+1})-E_h(n^k)
  &\leq h\sum_{j=1}^N
    \mu_j^{k+1}(n_j^{k+1}-n_j^k)\\
  &=-\tau h\sum_{j=1}^{N-1}
    n_{j+1/2}^k
    |(\delta_hp^{k+1})_{j+1/2}|^2 -\tau h\sum_{j=1}^N
    M_j^k|\mu_j^{k+1}|^2,
\end{aligned}
\]
which proves \eqref{eq:multi-discrete-energy-ineq}.
\end{proof}

\begin{proposition}[Species reaction step]
\label{prop:multi-species-reaction}
Under the hypotheses of
Proposition~\ref{prop:multi-scalar-step},
the reaction update
\eqref{eq:multi-post-reaction_2}
is well defined and satisfies
\begin{equation}\label{prop:n_species}
  \widehat{\boldsymbol n}_j^{k+1}\geq\boldsymbol0,
  \qquad
  \mathbf1^T\widehat{\boldsymbol n}_j^{k+1}
  =\widehat n_j^{k+1}
\end{equation}
in every cell.
\end{proposition}

\begin{proof}
By Proposition~\ref{prop:multi-scalar-step},
we have $0\leq p_j^{k+1}\leq p_H$
and $\widehat n_j^{k+1}\geq0$.
For $n_j^k>0$, the pressure bound gives
$R_{\ell i}(p_j^{k+1})\geq0$ for $\ell\ne i$.
Choose $c_j\geq0$ such that
$\mathbf R(p_j^{k+1})^T+c_j\mathbf I$
is entrywise nonnegative. Then
\[
  \exp\!\bigl(\tau\mathbf R(p_j^{k+1})^T\bigr)
  \geq e^{-c_j\tau}\mathbf I
\]
entrywise.
Since $\boldsymbol\phi_j^k$ is nonnegative
and has component sum one, the numerator in
\eqref{eq:multi-reaction-composition} is nonnegative,
and its denominator is at least $e^{-c_j\tau}>0$.
Thus the normalization is well defined and gives
$\widehat{\boldsymbol\phi}_j^{k+1}\in\Delta_m$,
where $\Delta_m$ is defined in \eqref{eq:simplex}.
Multiplication by $\widehat n_j^{k+1}\geq0$
in \eqref{eq:multi-post-reaction_2}
therefore yields \eqref{prop:n_species}.

For $n_j^k=0$, both $\widehat n_j^{k+1}$
and $\widehat{\boldsymbol n}_j^{k+1}$ vanish
by the zero-density case treated in Step~2.
Hence the conclusions hold in every cell.
\end{proof}

\begin{proposition}[Species transport step]
\label{prop:multi-species-remapping}
Under the hypotheses of
Proposition~\ref{prop:multi-scalar-step},
the transport system
\eqref{eq:multi-species-conservative-step}
is uniquely solvable on the active set
$\mathcal A^{k+1}$, with
$\boldsymbol\phi_j^{k+1}\in\Delta_m$
for every $j\in\mathcal A^{k+1}$.
The reconstructed species densities satisfy
\begin{equation}\label{eq:multi-species-properties}
  \boldsymbol n_j^{k+1}\geq\boldsymbol0,
  \qquad
  \mathbf1^T\boldsymbol n_j^{k+1}=n_j^{k+1}
\end{equation}
in every cell.
Moreover, the species fluxes satisfy
\begin{equation}\label{eq:multi-flux-compatibility}
  \mathbf1^T\boldsymbol F_{j+1/2}^{k+1}
  =F_{j+1/2}^{k+1}.
\end{equation}
\end{proposition}

\begin{proof}
Under the pressure ordering specified in Step~3,
the system \eqref{eq:multi-remapping-recursion}
is triangular, with positive diagonal entries
\(D_j^{k+1}\) on \(\mathcal A^{k+1}\).
It therefore uniquely determines
\(\boldsymbol\phi_j^{k+1}\) on the active set.

To prove that these compositions belong to \(\Delta_m\),
we combine the scalar balance \eqref{eq:multi-density-conservative-step}
with the definition of \(D_j^{k+1}\)
in \eqref{eq:multi-active-diagonal} to obtain
\begin{equation}\label{eq:multi-remapping-convex}
  D_j^{k+1}
  =
  \widehat n_j^{k+1}
  +\frac{\tau}{h}
  \left[
    (F_{j+1/2}^{k+1})^+
    -(F_{j-1/2}^{k+1})^-
  \right].
\end{equation}
We proceed by induction in decreasing pressure order.
On the right-hand side of
\eqref{eq:multi-remapping-recursion},
the coefficients of
\(\boldsymbol\phi_{j+1}^{k+1}\) and
\(\boldsymbol\phi_{j-1}^{k+1}\) are nonnegative.
A coefficient is nonzero only when the corresponding
neighboring cell has higher pressure than cell \(j\).
Thus every composition that contributes to these terms
has already been computed and belongs to \(\Delta_m\)
by the induction hypothesis.
For an active cell of maximal pressure,
both coefficients vanish, so no neighboring composition
is needed.
By \eqref{prop:n_species},
\(\widehat{\boldsymbol n}_j^{k+1}\geq0\).
Hence the right-hand side of
\eqref{eq:multi-remapping-recursion} is nonnegative.
Taking the component sum of this equation and using
\eqref{prop:n_species} gives
\[
  D_j^{k+1}\mathbf1^T\boldsymbol\phi_j^{k+1}
  =
  \widehat n_j^{k+1}
  +\frac{\tau}{h}
  \left[
    (F_{j+1/2}^{k+1})^+
    -(F_{j-1/2}^{k+1})^-
  \right]
  =D_j^{k+1},
\]
where the last equality follows from
the scalar identity for \(D_j^{k+1}\) above.
Division by \(D_j^{k+1}>0\) therefore yields
\(\boldsymbol\phi_j^{k+1}\geq0\) and
\(\mathbf1^T\boldsymbol\phi_j^{k+1}=1\).
This completes the induction.

The reconstruction
\eqref{eq:multi-species-reconstruction}
uniquely determines the species densities
and gives \eqref{eq:multi-species-properties}
on $\mathcal A^{k+1}$.
Outside the active set, both the total density
and the species densities vanish,
as established in Step~3.
Thus \eqref{eq:multi-species-properties}
holds in every cell.

Finally, every nonzero flux in
\eqref{eq:multi-upwind-species-flux}
uses the composition of an active upwind cell.
Since this composition has component sum one,
summing \eqref{eq:multi-upwind-species-flux}
over the species gives
\eqref{eq:multi-flux-compatibility}.
The same identity holds when the total flux is zero
and at the boundary interfaces.
\end{proof}

\begin{corollary}[Structure preservation of the full scheme]
\label{cor:multi-structure-preservation}
Suppose that
\[
  n_{i,j}^0\geq0,
  \qquad
  n_j^0=\sum_{i=1}^m n_{i,j}^0,
  \qquad
  0\leq p_j^0=(n_j^0)^\gamma\leq p_H.
\]
For any $\tau>0$, the fully discrete scheme
is uniquely solvable at every time step
and satisfies, for all $k\geq0$,
\begin{equation}\label{eq:multi-discrete-density-bounds}
  \begin{aligned}
    n_{i,j}^k \geq0,
    &\quad
    \mathbf1^T\boldsymbol n_j^k =n_j^k, \quad
    0\leq n_j^k \leq n_H,
    &\quad
    0\leq p_j^k \leq p_H.
  \end{aligned}
\end{equation}
The total mass satisfies
\begin{equation}\label{eq:multi-discrete-mass-monotonicity}
  0
  \leq h\sum_{j=1}^N n_j^k
  \leq h\sum_{j=1}^N n_j^{k+1}
  \leq n_H|\Omega|,
\end{equation}
and the energy inequality
\eqref{eq:multi-discrete-energy-ineq}
holds at every time step.
\end{corollary}

\begin{proof}
The result follows by induction from
Propositions~\ref{prop:multi-scalar-step},
\ref{prop:multi-species-reaction}, and
\ref{prop:multi-species-remapping}.
\end{proof}

\subsection{Second-order temporal extension}
\label{subsec:multi-second-order}

We construct a predictor--corrector extension of the first-order
scheme in Section~\ref{subsec:multi-full-discrete},
keeping the spatial discretization unchanged.
A first-order half step provides the midpoint prediction
for the scalar and species correctors.
If the trial update violates nonnegativity or the pressure bound,
we revert to the first-order scheme with the original time step.

\paragraph{Step 1: Midpoint predictor.}
Let $\mathcal S_\tau$ denote one complete step of the
first-order scheme with time step $\tau$.
We first apply the first-order scheme over half a time
step to obtain
\begin{equation}\label{eq:multi-order2-predictor}
  \boldsymbol n^P
  =\mathcal S_{\tau/2}(\boldsymbol n^k),
  \qquad
  n_j^P=\mathbf1^T\boldsymbol n_j^P.
\end{equation}
By Corollary~\ref{cor:multi-structure-preservation},
the predicted densities satisfy
$
  \boldsymbol n_j^P\geq\boldsymbol0
$ and
$
  0\leq n_j^P\leq n_H.
$

\paragraph{Step 2: Scalar corrector.}

We compute the trial total density
$n_j^{\mathrm{tr}}$ using the midpoint prediction
\eqref{eq:multi-order2-predictor}.
To derive a discrete energy dissipation law,
we introduce the averaged pressure
\begin{equation}\label{eq:multi-order2-dvd}
  \overline p_j
  =
  \int_0^1
  \bigl[(1-s)n_j^k+s n_j^{\mathrm{tr}}\bigr]^\gamma
  \,ds
  =
  \begin{cases}
    \displaystyle
    \frac{(n_j^{\mathrm{tr}})^{\gamma+1}
          -(n_j^k)^{\gamma+1}}
         {(\gamma+1)(n_j^{\mathrm{tr}}-n_j^k)},
       & n_j^{\mathrm{tr}}\ne n_j^k,\\[2ex]
    (n_j^k)^\gamma,
       & n_j^{\mathrm{tr}}=n_j^k.
  \end{cases}
\end{equation}
As in \eqref{eq:mu}, we set
$\overline\mu_j=\overline p_j-p_H$.
For the discrete energy $E_h$ defined in
\eqref{eq:multi-discrete-energy},
this choice gives the exact chain rule
\begin{equation}\label{eq:E_chain_rule}
  E_h(n^{\mathrm{tr}})-E_h(n^k)
  =
  h\sum_{j=1}^N
    \overline\mu_j(n_j^{\mathrm{tr}}-n_j^k).
\end{equation}

Using the predicted densities, we define
the transport coefficient and reaction mobility by
\begin{equation}\label{eq:multi-order2-mobilities}
  n_{j+1/2}^P
  =\frac{n_j^P+n_{j+1}^P}{2},
  \qquad
  M_j^P
  =-\sum_{\ell=1}^m
    n_{\ell,j}^P\sigma_\ell(\overline p_j).
\end{equation}
The expression for $M_j^P$ follows
\eqref{eq:multi-reaction-mobility},
with the pressure evaluated at $\overline p_j$.

The trial total density satisfies
\begin{equation}\label{eq:multi-order2-scalar}
  \frac{n_j^{\mathrm{tr}}-n_j^k}{\tau}
  =
  (d_hF^{\mathrm{tr}})_j
  -M_j^P\overline\mu_j,
  \qquad
  F_{j+1/2}^{\mathrm{tr}}
  =
  n_{j+1/2}^P
  (\delta_h\overline p)_{j+1/2},
\end{equation}
with
$F_{1/2}^{\mathrm{tr}}
=F_{N+1/2}^{\mathrm{tr}}=0$.
During the solve, the predicted densities
and $n_{j+1/2}^P$ remain fixed,
while $M_j^P$ depends on the trial density through
$\overline p_j$.
Thus \eqref{eq:multi-order2-scalar}
is a nonlinear system for the total density alone
and can be solved by Newton's method.

\paragraph{Step 3: Species corrector.}

We compute the trial species densities using
the midpoint prediction in Step~1 and the total flux
obtained in Step~2.
Unlike the first-order scheme, which treats reaction
and transport successively, we combine them in a
single corrector step to make direct use of
the midpoint prediction for both contributions.

For transport, we retain the upwind formula
\eqref{eq:multi-upwind-species-flux}
and evaluate it using the predicted composition.
For $n_j^P>0$, set
$\boldsymbol\phi_j^P=\boldsymbol n_j^P/n_j^P$.
The predictor \eqref{eq:multi-order2-predictor}
gives $\boldsymbol\phi_j^P\in\Delta_m$, where $\Delta_m$ is defined in \eqref{eq:simplex}.
These compositions are known from Step~1.
Whenever the required compositions are defined,
we set
\begin{equation}\label{eq:multi-order2-species-flux}
  \boldsymbol F_{j+1/2}^{\mathrm{tr}}
  =
  (F_{j+1/2}^{\mathrm{tr}})^-\boldsymbol\phi_j^P
  +(F_{j+1/2}^{\mathrm{tr}})^+\boldsymbol\phi_{j+1}^P.
\end{equation}
At the boundaries, we set
$\boldsymbol F_{1/2}^{\mathrm{tr}}
=\boldsymbol F_{N+1/2}^{\mathrm{tr}}
=\boldsymbol0$.
No composition is needed when
$F_{j+1/2}^{\mathrm{tr}}=0$.
However, if $F_{j+1/2}^{\mathrm{tr}}\neq0$,
the composition required by
\eqref{eq:multi-order2-species-flux}
must be defined, which is a condition that will be checked
in Step~4.

The reaction rates are evaluated at the predicted
densities and averaged pressure
\begin{equation}\label{eq:multi-order2-reaction}
  \boldsymbol q_j^P
  =
  \mathbf R(\overline p_j)^T\boldsymbol n_j^P,
\end{equation}
where $\mathbf R$ is the matrix of reaction
coefficients in \eqref{eq:intro-model}.
This term is well defined when $n_j^P=0$,
in which case $\boldsymbol q_j^P=\boldsymbol0$.

We combine the transport and reaction contributions
to compute the trial species densities via
\begin{equation}\label{eq:multi-order2-species}
  \frac{\boldsymbol n_j^{\mathrm{tr}}
        -\boldsymbol n_j^k}{\tau}
  =
  (d_h\boldsymbol F^{\mathrm{tr}})_j
  +\boldsymbol q_j^P.
\end{equation}
This update is explicit once the scalar corrector
has been solved.

\paragraph{Step 4: Admissibility and fallback.}

If the trial densities are obtained, we check their admissibility.
To ensure that the species flux
\eqref{eq:multi-order2-species-flux}
is well defined, we require $n_j^P>0$
when $F_{j+1/2}^{\mathrm{tr}}<0$,
and $n_{j+1}^P>0$
when $F_{j+1/2}^{\mathrm{tr}}>0$,
at every interior interface.
The trial densities must also satisfy
\[
  0\leq n_j^{\mathrm{tr}}\leq n_H,
  \qquad
  \boldsymbol n_j^{\mathrm{tr}}\geq\boldsymbol0
  \quad\text{for every }j.
\]

If these conditions hold, we accept the trial.
Otherwise, we discard it and apply the complete
first-order scheme from $\boldsymbol n^k$
with the original time step $\tau$.
To be more specific, we choose
\begin{equation}\label{eq:multi-order2-accept}
  \boldsymbol n^{k+1}
  =
  \begin{cases}
    \boldsymbol n^{\mathrm{tr}},
      &\text{if the trial is admissible},\\
    \mathcal S_\tau(\boldsymbol n^k),
      &\text{otherwise}.
  \end{cases}
\end{equation}

\begin{proposition}[Structure preservation]
\label{prop:multi-order2-structure}
Under the hypotheses of
Proposition~\ref{prop:multi-scalar-step},
the update \eqref{eq:multi-order2-accept} satisfies,
for any $\tau>0$,
\[
  \boldsymbol n_j^{k+1}\geq\boldsymbol0,
  \qquad
  \mathbf1^T\boldsymbol n_j^{k+1}
  =n_j^{k+1}\in[0,n_H],
  \qquad
  E_h(n^{k+1})\leq E_h(n^k),
\]
where $E_h$ is defined in
\eqref{eq:multi-discrete-energy}.
The total mass satisfies
\[
  h\sum_{j=1}^N n_j^k
  \leq h\sum_{j=1}^N n_j^{k+1}
  \leq n_H|\Omega|.
\]
For an accepted trial, the reaction rates
and species fluxes satisfy
\begin{equation}\label{eq:multi-order2-compatibility}
  \mathbf1^T\boldsymbol q_j^P
  =-M_j^P\overline\mu_j,
  \qquad
  \mathbf1^T\boldsymbol F_{j+1/2}^{\mathrm{tr}}
  =F_{j+1/2}^{\mathrm{tr}},
\end{equation}
and the discrete energy satisfies
\begin{equation}\label{eq:multi-order2-energy}
  E_h(n^{k+1})-E_h(n^k)
  =
  -\tau h\sum_{j=1}^{N-1}
    n_{j+1/2}^P
    |(\delta_h\overline p)_{j+1/2}|^2
  -\tau h\sum_{j=1}^N
    M_j^P|\overline\mu_j|^2
  \leq0.
\end{equation}
\end{proposition}

\begin{proof}
If the trial is rejected,
\eqref{eq:multi-order2-accept} selects
the first-order update, and the conclusions follow
from Corollary~\ref{cor:multi-structure-preservation}.
We therefore consider an accepted trial.

The admissibility conditions in Step~4 give
$\boldsymbol n_j^{\mathrm{tr}}\geq\boldsymbol0$
and $n_j^{\mathrm{tr}}\in[0,n_H]$.
Since $n_j^k\in[0,n_H]$,
the averaged pressure
\eqref{eq:multi-order2-dvd} satisfies
$0\leq\overline p_j\leq p_H$,
so $\overline\mu_j\leq0$.
The predictor \eqref{eq:multi-order2-predictor}
is nonnegative by
Corollary~\ref{cor:multi-structure-preservation},
and $\sigma_\ell(\overline p_j)\leq0$  follows from \eqref{eq:intro-common-ph}
and \eqref{eq:continuous-secant-coefficient}.
Thus \eqref{eq:multi-order2-mobilities} gives
$n_{j+1/2}^P\geq0$ and $M_j^P\geq0$.

We first verify
\eqref{eq:multi-order2-compatibility}.
The reaction rate
\eqref{eq:multi-order2-reaction}
and the definition \eqref{eq:intro-total-equation} give
\[
  \mathbf1^T\boldsymbol q_j^P
  =
  \sum_{\ell=1}^m
    n_{\ell,j}^P H_\ell(\overline p_j)
  =
  (\overline p_j-p_H)
  \sum_{\ell=1}^m
    n_{\ell,j}^P\sigma_\ell(\overline p_j)
  =
  -M_j^P\overline\mu_j,
\]
where the last two equalities follow from
\eqref{eq:continuous-secant-coefficient},
\eqref{eq:multi-order2-mobilities},
and $\overline\mu_j=\overline p_j-p_H$
as in \eqref{eq:mu}.
For each nonzero transport flux, Step~4 ensures that
the required predicted composition is defined and
has component sum one.
Summing the species fluxes therefore gives
\(\mathbf1^T\boldsymbol F_{j+1/2}^{\mathrm{tr}}
=F_{j+1/2}^{\mathrm{tr}}\).
This identity also holds at zero-flux interfaces
and at the boundaries.

Summing \eqref{eq:multi-order2-species}
over the species and subtracting
\eqref{eq:multi-order2-scalar}, we obtain,
by \eqref{eq:multi-order2-compatibility},
\[
  \mathbf1^T\boldsymbol n_j^{\mathrm{tr}}
  -n_j^{\mathrm{tr}}
  =
  \mathbf1^T\boldsymbol n_j^k-n_j^k
  =0.
\]
This proves exact reconstruction of the total density.

Summing \eqref{eq:multi-order2-scalar}
over all cells gives the total-mass monotonicity,
since the boundary fluxes vanish and
$-M_j^P\overline\mu_j\geq0$.
The upper bound follows from
$n_j^{\mathrm{tr}}\leq n_H$.

Finally, multiplying
\eqref{eq:multi-order2-scalar}
by $\tau h\overline\mu_j$,
summing over all cells, and applying discrete
summation by parts gives
\eqref{eq:multi-order2-energy},
using $\delta_h\overline\mu=\delta_h\overline p$,
the discrete chain rule \eqref{eq:E_chain_rule},
and $n^{k+1}=n^{\mathrm{tr}}$.
The nonnegativity of $n_{j+1/2}^P$ and $M_j^P$
then yields energy dissipation.
\end{proof}

\begin{remark}
For smooth solutions, the predictor and the averaged pressure
in \eqref{eq:multi-order2-dvd} approximate their midpoint values
to $O(\tau^2)$, yielding an $O(\tau^3)$ one-step error
and hence formal second-order accuracy in time.
This order is not guaranteed when the method reverts
to the first-order scheme.
\end{remark}

\section{Numerical experiments}
\label{sec:numerical}

We examine the accuracy, computational efficiency and structure-preserving properties of the proposed schemes. We also investigate the effects of species composition and pressure-dependent switching on tumor growth, as well as the behavior in the stiff-pressure limit, through one- and two-dimensional simulations.

Unless otherwise stated, we set \(p_H=1\). When \(n^0\) and
\(\boldsymbol\phi^0\) are prescribed separately, we initialize the
species densities from \(n_i^0=n^0\phi_i^0\).
The nonlinear systems \eqref{eq:multi-discrete-density-step} and \eqref{eq:multi-order2-scalar} are solved by a damped Newton method with backtracking,
with the stopping tolerance for the scaled residuals set to $10^{-12}$.

\subsection{Accuracy and efficiency}
\label{subsec:accuracy-structure}

We examine the temporal accuracy of the first-order scheme and its second-order predictor--corrector
extension using both the smooth positive initial data and the near-vacuum
initial data with steep interfaces. We then compare their
computational efficiency.

\medskip
\noindent\textit{Smooth positive initial data.}
We take $\Omega=(-1,1)$ and
\begin{equation}\label{eq:case_s}
 n^0(x)=0.55+0.15\cos(\pi x),\qquad
 \phi_1^0(x)=0.5+0.2\cos(2\pi x),\qquad
 \phi_2^0(x)=1-\phi_1^0(x).
\end{equation}
The reaction matrix is
\[
\mathbf R(p)=
\begin{pmatrix}
 H_1(p)-\kappa_{12}(p)&\kappa_{12}(p)\\
 \kappa_{21}(p)&H_2(p)-\kappa_{21}(p)
\end{pmatrix},
\]
where
\[
\begin{aligned}
H_1(p)=2(1-p)(1+0.2p),
\quad  H_2(p)=(1-p)(1+0.5p),
\quad \kappa_{12}(p)=0.2p,
\quad \kappa_{21}(p)=0.1(1-p)^2.
\end{aligned}
\]

\medskip
\noindent\textit{Near-vacuum initial data.}
We take $\Omega=(-10,10)$ and
\begin{equation}\label{eq:case_a}
\begin{aligned}
 n^0(x)&=0.49\left[
 \tanh\left(\frac{x+2}{0.05}\right)
 -\tanh\left(\frac{x-2}{0.05}\right)\right],\\
 \phi_1^0(x)&=\frac12\left[
 1-\tanh\left(\frac{x-0.25}{0.03}\right)\right],
 \qquad \quad  \phi_2^0(x)=1-\phi_1^0(x),
\end{aligned}
\end{equation}
with
$\mathbf R(p)=\operatorname{diag}(10(1-p),2(1-p))$.

\paragraph{Temporal convergence.}
We take $\gamma=20$, $N=1000$, and compute till  $T=0.2$.
The reference solutions are computed using the
predictor--corrector scheme with
$\tau_{\rm ref}=5\times10^{-6}$.
We measure the relative species-density error by

\[
\mathrm{err}_t=
\frac{\|\boldsymbol n_h(T)-\boldsymbol n_{\rm ref}(T)\|_{1,h}}
{\|\boldsymbol n_{\rm ref}(T)\|_{1,h}},
\qquad
\|\boldsymbol v\|_{1,h}
=h\sum_{i=1}^m\sum_{j=1}^N |v_{i,j}|.
\]

\begin{figure}[htb]
\centering
\begin{minipage}[t]{0.47\linewidth}
  \centering
  \includegraphics[width=\linewidth,
    trim=0 0.4in 0 0.4in,clip]{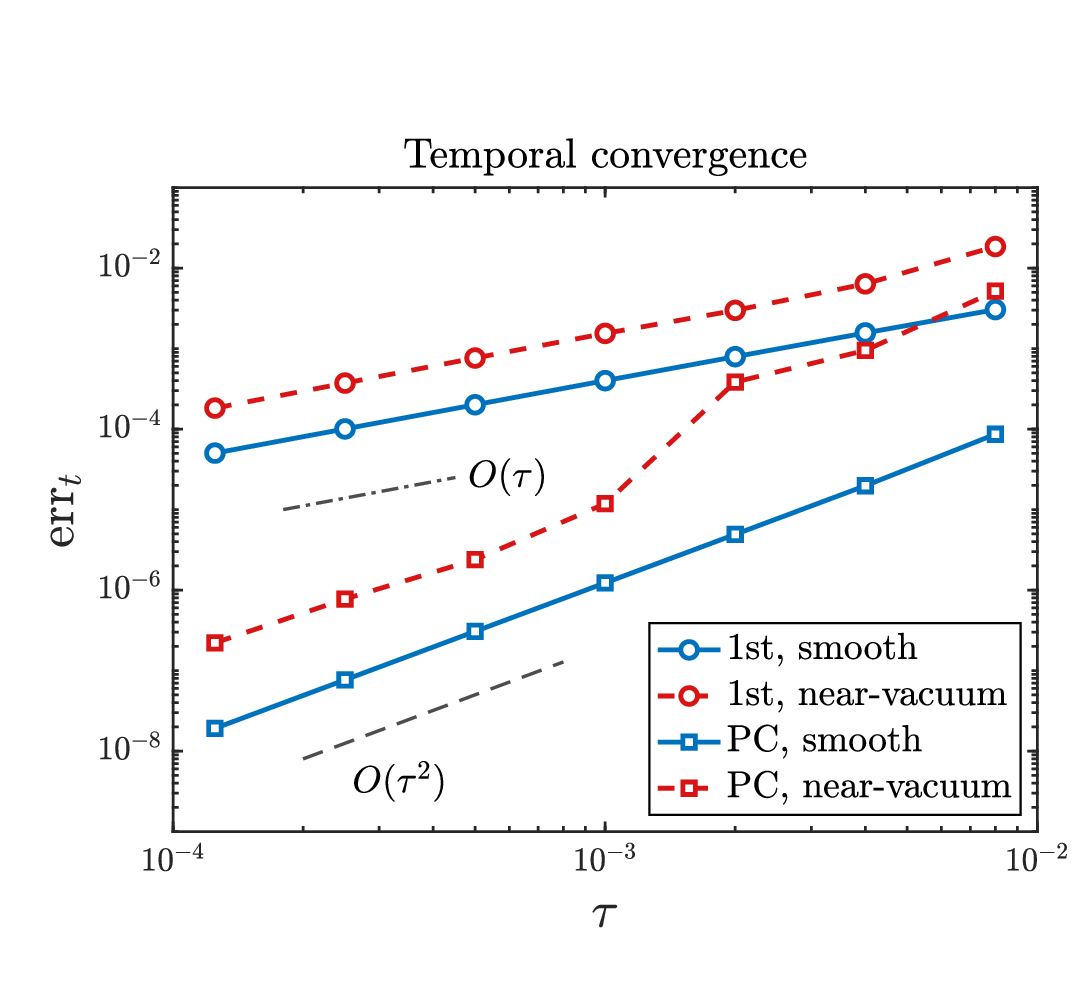}
\end{minipage}\hfill
\begin{minipage}[t]{0.47\linewidth}
  \centering
  \includegraphics[width=\linewidth,
    trim=0 0.4in 0 0.4in,clip]{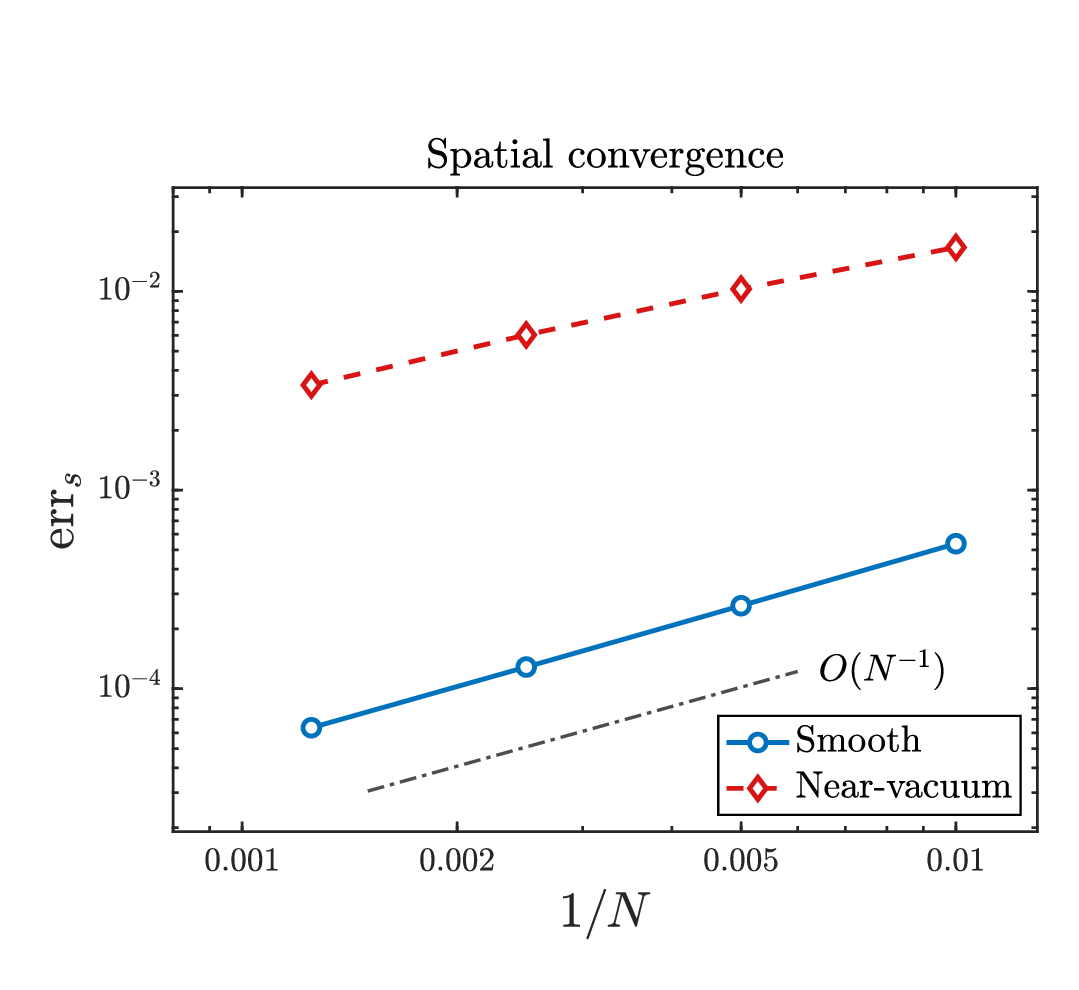}
\end{minipage}
\caption{Temporal (left) and spatial (right) convergence
for the smooth and near-vacuum tests.
Here `1st' denotes the first-order scheme and `PC' denotes the second-order predictor--corrector scheme.
Grey lines indicate reference slopes.}
\label{fig:multi-convergence}
\end{figure}

The left panel of Figure~\ref{fig:multi-convergence} shows
approximately first- and second-order temporal convergence
for the first-order and predictor--corrector schemes,
respectively, in both tests,
although its observed order is slightly below two
for the near-vacuum test.

\paragraph{Spatial convergence.}
We use the first-order scheme with $\gamma=20$,
$T=0.2$, and $\tau=10^{-5}$ on successively refined grids.
We measure the relative spatial error by
\[
\mathrm{err}_s=
\frac{\|\boldsymbol n_h(T)-\boldsymbol n_{h/2}(T)\|_{1,h}}
{\|\boldsymbol n_{h/2}(T)\|_{1,h}}.
\]

The right panel of Figure~\ref{fig:multi-convergence}  shows
approximately first-order spatial convergence for both
the smooth and near-vacuum tests.

\paragraph{Computational efficiency.}
For the smooth test \eqref{eq:case_s}, we compare the computation times of
the first-order and predictor--corrector schemes with
$N=1000$ and $T=0.2$.
For $\gamma=3$ and $20$, we fix the time step of the
predictor--corrector scheme to be $0.01$ and adjust that
of the first-order scheme to achieve similar temporal errors.
The reference solutions are computed using the
predictor--corrector scheme with
$\tau_{\mathrm{ref}}=1/25600$.
Table~\ref{tab:temporal-efficiency} presents the temporal
errors and computation times.
The predictor--corrector scheme is more efficient
in both cases, as its higher temporal accuracy allows
much larger time steps.

\begin{table}[htbp]
\centering
\small
\caption{Comparison of computation times for the smooth test \eqref{eq:case_s}
at similar temporal errors.
Speedup is the first-order computation time divided by
the predictor--corrector computation time.}
\label{tab:temporal-efficiency}
\vspace{2mm}
\setlength{\tabcolsep}{15pt}
\renewcommand{\arraystretch}{1.4}
\begin{tabular}{@{}c|c|ccc|c@{}}
\hline
$\gamma$ & Scheme & $\tau$ & $\mathrm{err}_t$
& Time (s) & Speedup\\
\hline
\multirow{2}{*}{$3$}
& first-order & $3.90\mathrm{E}{-5}$ & $1.38\mathrm{E}{-5}$
& $23.9$ & \multirow{2}{*}{$118$}\\
\cline{2-5}
& predictor--corrector  & $1.00\mathrm{E}{-2}$ & $1.37\mathrm{E}{-5}$
& $0.202$ &\\
\hline
\multirow{2}{*}{$20$}
& first-order & $3.51\mathrm{E}{-4}$ & $1.41\mathrm{E}{-4}$
& $10.8$ & \multirow{2}{*}{$8.77$}\\
\cline{2-5}
& predictor--corrector  & $1.00\mathrm{E}{-2}$ & $1.39\mathrm{E}{-4}$
& $1.23$ &\\
\hline
\end{tabular}
\end{table}

\subsection{Structure preservation}
\label{subsec:structure-preservation}

We examine the structural properties of the scheme
and the evolution of the total density, common pressure
and species densities in two- and three-species tests
without switching. Both tests use the first-order scheme.
We first consider the two-species near-vacuum test
\eqref{eq:case_a} with $\Omega=(-10,10)$, $\gamma=20$,
$N=2000$, $\tau=10^{-4}$ and $T=2$.

\paragraph{Density, pressure, energy and mass.}

\begin{figure}[htbp]
  \centering
  \includegraphics[width=0.99\linewidth,  trim=0 .8cm 0 0cm,
    clip]{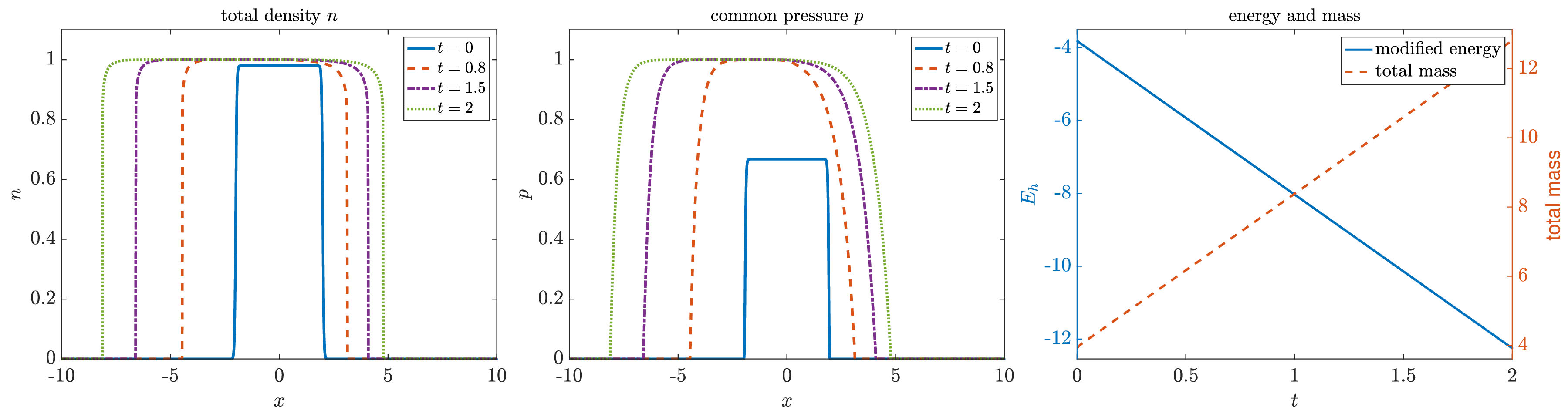}
  \caption{
    Evolution of the total density $n$ (left) and common
pressure $p=n^\gamma$ (middle), together with the
modified energy $E_h$ and total mass (right), for the
two-species near-vacuum test \eqref{eq:case_a} with $\gamma=20$.
  }
  \label{fig:structure-density-pressure}
\end{figure}

Figure~\ref{fig:structure-density-pressure} shows the evolution
of the total density and common pressure, together with the
modified energy and total mass.
The density profile expands outward, with the left front
advancing faster than the right, consistent with the higher
growth rate of the species on the left.
Both the total density and common pressure approach their
upper bounds $n_H=p_H=1$ in the interior.
The modified energy decreases, while the total mass increases.

\paragraph{Two- and three-species comparison.}

We construct a three-species test by splitting the
faster-growing species of the two-species model into
two species with the same growth rate. The reaction
matrix is
\[
  \mathbf R(p)=\operatorname{diag}
  \bigl(10(1-p),\,10(1-p),\,2(1-p)\bigr).
\]
We use the same initial total density and numerical
parameters as in the two-species test.
The initial species fractions are
\[
  \phi_1^0=f(1-s),\qquad
  \phi_2^0=fs,\qquad
  \phi_3^0=1-f,
\]
where $f(x)=\phi_1^0(x)$ from \eqref{eq:case_a} and
\(
  s(x)=\frac12\left[
    1+\tanh\left(\frac{x+0.9}{0.05}\right)
  \right].
\)
The aggregated densities $n_1+n_2$ and $n_3$ satisfy
the original two-species equations with the same
initial data. Hence, the two models have the same
total density and common pressure at the continuous
level.

\begin{figure}[htbp]
  \centering
  \includegraphics[width=0.99\linewidth,  trim=0 .8cm 0 0cm,
    clip]{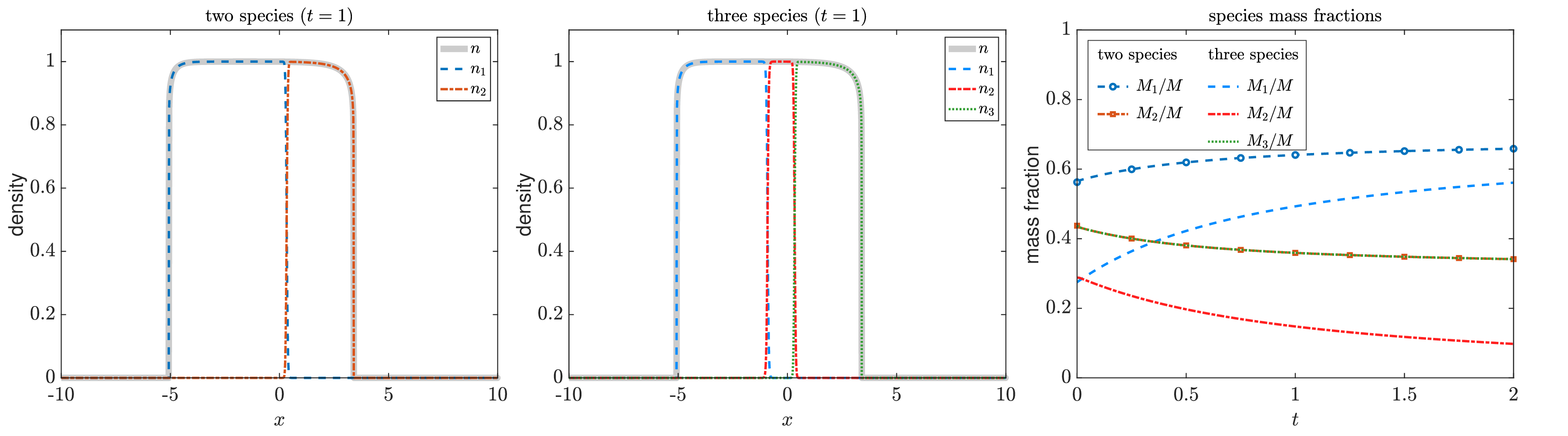}
  \caption{
    Species density profiles at $t=1$ for the
    two-species (left) and three-species (middle)
    tests, and the corresponding species mass
    fractions (right).
  }
  \label{fig:structure-species}
\end{figure}

Figure~\ref{fig:structure-species} shows the species
densities at $t=1$ and the evolution of their mass fractions,
defined by
\[
  \frac{M_i(t^k)}{M(t^k)}
  =\frac{\sum_{j=1}^N n_{i,j}^k}
         {\sum_{j=1}^N n_j^k}.
\]
In both tests, the mass fraction of the leftmost species
increases, while that of the rightmost species decreases.
In the three-species test, the intermediate species also
loses mass fraction despite having the same growth
coefficient as the leftmost species.
The intermediate species lies mainly in the nearly
saturated interior, where growth is suppressed, whereas
the leftmost species continues to grow near the expanding
front.
The mass fractions of the rightmost species are nearly
identical in the two tests, as expected from the
aggregation relation.

In both tests, the species densities remain nonnegative,
the common pressure satisfies $p\le p_H$ and the
individual species masses are nondecreasing.

\subsection{Pressure-induced phenotype switching}
\label{subsec:num-transition}

We examine how pressure-dependent switching changes the spatial
distribution of species. To separate this effect from differences
in growth rates, we consider two phenotypes with the same net
growth rate and a spatially uniform initial composition.
We use the first-order scheme for this test.

We take $\Omega=(-10,10)$, $\gamma=20$, $N=2000$,
$\tau=10^{-4}$ and $T=2$.
The reaction matrix is
\[
  \mathbf R(p)
  =
  \begin{pmatrix}
    6(1-p)-\alpha p & \alpha p\\
    \beta(1-p) & 6(1-p)-\beta(1-p)
  \end{pmatrix}.
\]
Here, $\alpha=1$ and $\beta=0.25$.
The switching rates from type 1 to type 2 and from type 2
to type 1 are $\alpha p$ and $\beta(1-p)$, respectively.
Thus higher pressure promotes switching toward type 2,
whereas lower pressure favors the reverse transition.

We use the initial total density $n^0$ in
\eqref{eq:case_a}, with
\(
  \phi_1^0(x)=\phi_2^0(x)=\frac12.
\)
Since both species have the net growth rate $6(1-p)$,
the total-density equation is independent of the composition.
Without switching, the fractions would remain equal to
$1/2$ wherever the total density is positive.
Adding switching therefore changes the composition while
leaving the total-density and pressure evolution unchanged.

For a fixed pressure, the equilibrium fractions of the
switching dynamics are
\[
  \phi_1^*(p)
  =
  \frac{\beta(1-p)}
       {\alpha p+\beta(1-p)},
  \qquad
  \phi_2^*(p)
  =
  \frac{\alpha p}
       {\alpha p+\beta(1-p)}.
\]
The equilibrium fraction $\phi_2^*(p)$  increases with pressure. The actual composition also depends on transport and the evolving pressure field.

\begin{figure}[htbp]
  \centering
  \includegraphics[width=0.6\textwidth,  trim=0 .7cm 0 0cm,
    clip]{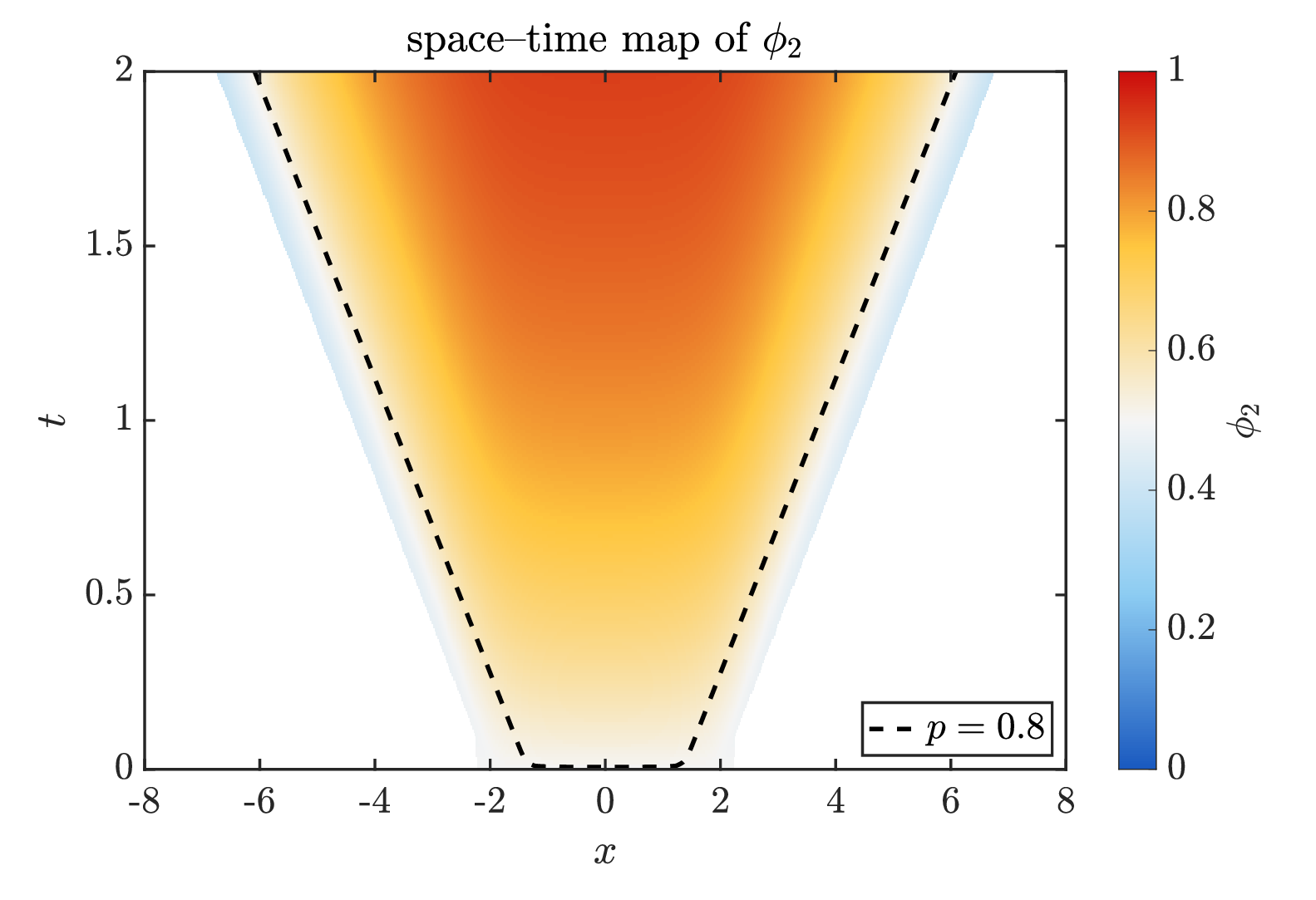}
  \caption{Space--time evolution of the type-2 fraction
  $\phi_2$ under pressure-dependent switching.
  The composition is shown only where $n>10^{-4}$,
  and the black curve denotes the pressure contour $p=0.8$.}
  \label{fig:example2-pressure-switching}
\end{figure}

Figure~\ref{fig:example2-pressure-switching} shows that the fraction
of type 2 is higher in the high-pressure interior than
near the moving front, consistent with the pressure
dependence of $\phi_2^*(p)$. Although both species have
the same net growth rate, pressure-dependent switching
causes the initially uniform composition to become
spatially nonuniform.

\subsection{Stiff-pressure limit as $\gamma\to\infty$}
\label{subsec:num-stiff-pressure}

We next examine the approach to the Hele--Shaw limit as
$\gamma$ increases while the spatial grid and time step
remain fixed.
We compare the total density and pressure with an explicit
limiting solution, and use the species masses to assess
the limiting relation for the internal species interface.
We use the first-order scheme for this test.

We consider the diagonal reaction matrix
\[
	 \mathbf R(p) = \operatorname{diag}
  \bigl(1-p,\,1-p\bigr).
\]
There is no switching, and both species have the net growth
rate $1-p$.
The total-density equation is therefore independent of
the composition.
For this scalar problem, the limiting outer boundary and
pressure are known explicitly
\cite{HuangRuan2026,
LiuTangWangZhou2018,PerthameQuirosVazquez2014}.
The internal species interface is transported by the
limiting Darcy velocity, consistently with the
two-population incompressible limit
\cite{BubbaPerthamePoucholSchmidtchen2020}.

We take $\Omega=(-5,5)$, $N=6400$, $\tau=10^{-4}$
and $T=1$. The initial data are
\begin{equation}\label{eq:initial-stiff}
  p^0(x)=\left(1-\frac{\cosh x}{\cosh R_0}\right)_+,
  \qquad
  n_\gamma^0(x)=[p^0(x)]^{1/\gamma},
\end{equation}
and
\[
  n_{1,\gamma}^0=n_\gamma^0\chi_{(-\infty,X_0)},
  \qquad
  n_{2,\gamma}^0=n_\gamma^0\chi_{[X_0,\infty)}.
\]
Here $R_0=1$ and $X_0=0.25$ are the initial tumor radius
and species interface, respectively. As $\gamma$ varies,
the initial pressure and interface remain fixed, while
the initial density changes.

The limiting outer boundary, total density and pressure are
\[
  R(t)
  =
  \operatorname{arsinh}\!\bigl(e^t\sinh R_0\bigr),\qquad
  n_\infty(x,t)
  =
  \chi_{[-R(t),R(t)]}(x), \qquad
  p_\infty(x,t)
  =
  \Big(
    1-\frac{\cosh x}{\cosh R(t)}
  \Big)_+.
\]
The internal interface satisfies
$\dot X=-\partial_xp_\infty(X(t),t)$ and $X(0)=X_0$,
which gives
\[
  X(t)
  =
  2\operatorname{artanh}\!\Big[
    \tanh\Big(\frac{X_0}{2}\Big)
    \frac{\tanh(R(t)/2)}{\tanh(R_0/2)}
  \Big].
\]
The limiting species densities are therefore
\[
  n_{1,\infty}(x,t)
  =
  \chi_{[-R(t),X(t)]}(x),
  \qquad
  n_{2,\infty}(x,t)
  =
  \chi_{[X(t),R(t)]}(x).
\]

\paragraph{Comparison with the limiting solution.}

\begin{figure}[htb]
  \centering
  \includegraphics[width=0.95\textwidth,  trim=0 .8cm 0 0cm,
    clip]{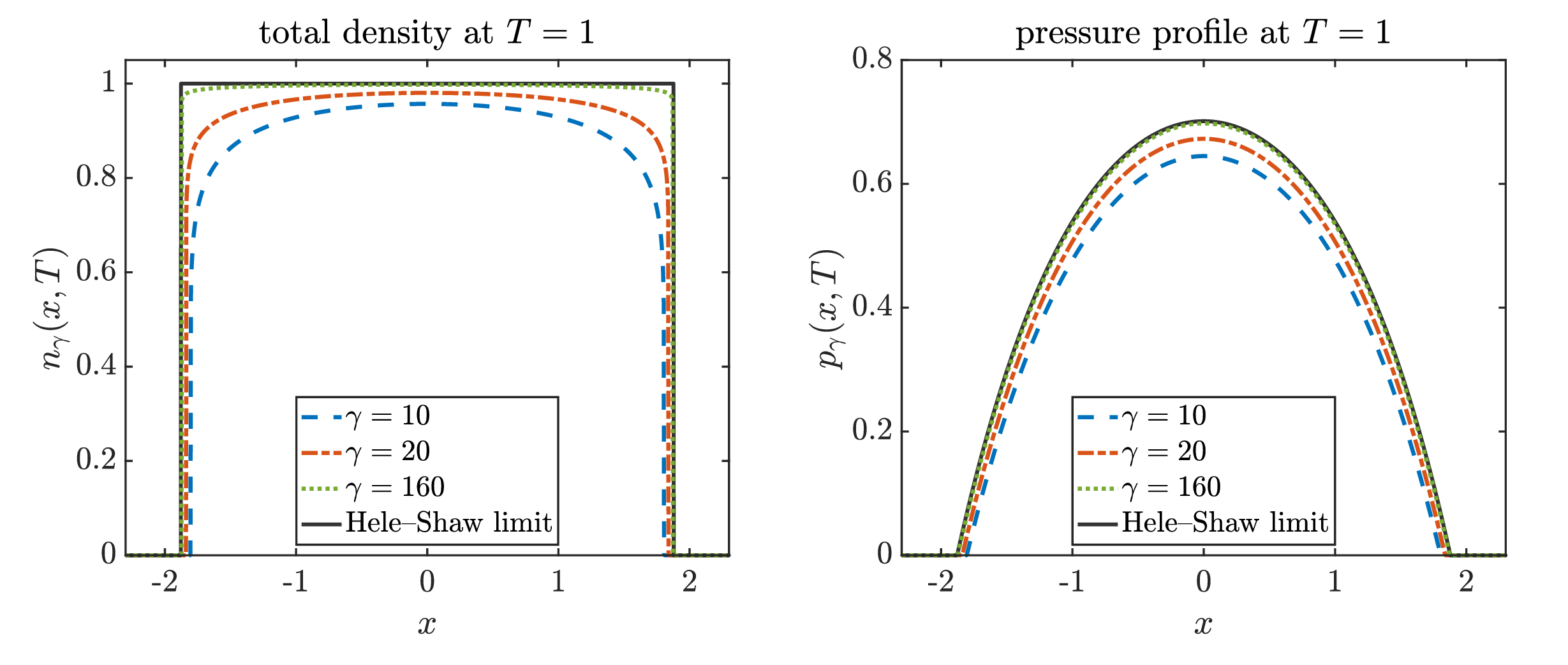}
\caption{
Total density (left) and pressure (right) at $T=1$
for $\gamma=10,20,160$.
The black solid curves denote the corresponding
Hele--Shaw limits.
All computations use $N=6400$ and $\tau=10^{-4}$.
}
  \label{fig:stiff-pressure-profiles}
\end{figure}

Figure~\ref{fig:stiff-pressure-profiles} shows that,
as $\gamma$ increases, the total density approaches
the characteristic function of the occupied region,
and the pressure approaches the Hele--Shaw profile.

For each $\gamma$, let $n_{\gamma,h}$ and $p_{\gamma,h}$
denote the numerical total density and pressure.
We measure the final-time density and pressure errors by
\begin{equation}
\label{eq:stiff-field-errors}
  e_n^\gamma
  =
  \|n_{\gamma,h}(T)-n_\infty(\cdot,T)\|_{1,h},
  \qquad
  e_p^\gamma
  =
  \|p_{\gamma,h}(T)-p_\infty(\cdot,T)\|_{1,h}.
\end{equation}
The limiting boundary and species interface also satisfy
the mass relations
\[
  R(T)
  =
  \frac12\int_\Omega n_\infty(x,T)\,\mathrm dx,
  \qquad
  X(T)
  =
  \frac12\int_\Omega
  \bigl(n_{1,\infty}(x,T)-n_{2,\infty}(x,T)\bigr)
  \,\mathrm dx.
\]
These relations motivate the mass-based estimates
\begin{equation}
\label{eq:stiff-interface-approximations}
  R_{\gamma,h}
  =
  \frac h2\sum_j n_{\gamma,j}(T),
  \qquad
  X_{\gamma,h}
  =
  \frac h2\sum_j
  \bigl(n_{1,\gamma,j}(T)-n_{2,\gamma,j}(T)\bigr),
\end{equation}
with errors
\begin{equation}
\label{eq:stiff-interface-errors}
  e_R^\gamma
  =
  |R_{\gamma,h}-R(T)|,
  \qquad
  e_X^\gamma
  =
  |X_{\gamma,h}-X(T)|.
\end{equation}
For finite $\gamma$, these quantities are obtained from
mass integrals and do not directly measure the geometric
locations of the numerical boundary and species interface.

Table~\ref{tab:stiff-pressure-fields} shows that all four
errors decrease approximately in proportion to $1/\gamma$
over the tested range.
The decay of $e_X^\gamma$ indicates that the numerical
species masses approach the relation defining the limiting
internal interface.
The observed orders for the density and mass-based estimates
decrease slightly at the largest values of $\gamma$.
Since the grid and time step remain fixed, these comparisons
include both spatial and temporal discretization errors.

\begin{table}[htb]
  \centering
  \small
  \caption{Density and pressure errors, together with errors
  in the mass-based boundary and interface estimates.
  Orders are measured with respect to $1/\gamma$.}
  \label{tab:stiff-pressure-fields}
  \vspace{3mm}
  \setlength{\tabcolsep}{6pt}
  \renewcommand{\arraystretch}{1.20}
  \begin{tabular}{c ccccccc}
    \hline
    $\gamma$
    & $10$ & $20$ & $40$ & $80$
    & $160$ & $320$ & $640$\\
    \hline
    $e_n^\gamma$
    & $4.76\mathrm{E}{-1}$
    & $2.45\mathrm{E}{-1}$
    & $1.25\mathrm{E}{-1}$
    & $6.34\mathrm{E}{-2}$
    & $3.23\mathrm{E}{-2}$
    & $1.67\mathrm{E}{-2}$
    & $8.83\mathrm{E}{-3}$\\
    order
    & -- & $0.956$ & $0.973$ & $0.978$
    & $0.973$ & $0.955$ & $0.917$\\
    $e_p^\gamma$
    & $2.27\mathrm{E}{-1}$
    & $1.19\mathrm{E}{-1}$
    & $6.12\mathrm{E}{-2}$
    & $3.10\mathrm{E}{-2}$
    & $1.55\mathrm{E}{-2}$
    & $8.02\mathrm{E}{-3}$
    & $4.04\mathrm{E}{-3}$\\
    order
    & -- & $0.928$ & $0.964$ & $0.980$
    & $0.998$ & $0.953$ & $0.988$\\
    \hline
    $e_R^\gamma$
    & $2.38\mathrm{E}{-1}$
    & $1.23\mathrm{E}{-1}$
    & $6.26\mathrm{E}{-2}$
    & $3.18\mathrm{E}{-2}$
    & $1.63\mathrm{E}{-2}$
    & $8.44\mathrm{E}{-3}$
    & $4.52\mathrm{E}{-3}$\\
    order
    & -- & $0.955$ & $0.972$ & $0.976$
    & $0.968$ & $0.946$ & $0.901$\\
    $e_X^\gamma$
    & $2.31\mathrm{E}{-2}$
    & $1.10\mathrm{E}{-2}$
    & $5.37\mathrm{E}{-3}$
    & $2.67\mathrm{E}{-3}$
    & $1.35\mathrm{E}{-3}$
    & $7.03\mathrm{E}{-4}$
    & $3.79\mathrm{E}{-4}$\\
    order
    & -- & $1.07$ & $1.03$ & $1.01$
    & $0.981$ & $0.946$ & $0.890$\\
    \hline
  \end{tabular}
\end{table}

\paragraph{Limiting structural relations.}
We also examine the complementarity relation and
the interior pressure equation
\[
  p_\infty(1-n_\infty)=0 \quad \text{in }\Omega,
  \qquad
  \partial_{xx}p_\infty+1-p_\infty=0
  \quad \text{in }\{p_\infty>0\}.
\]
We measure the defects in these two relations by
\[
  \mathcal C_\gamma
  =
  h\sum_j
  \left|
    p_{\gamma,j}(T)
    \bigl(1-n_{\gamma,j}(T)\bigr)
  \right|,\qquad
  \mathcal I_\gamma
  =
  h\sum_{j:\,p_\infty(x_j,T)\geq0.2}
  \left|
    (d_h\delta_h p_{\gamma,h}(T))_j
    +1-p_{\gamma,j}(T)
  \right|.
\]
The restriction in $\mathcal I_\gamma$ excludes a
neighborhood of the free boundary, where the limiting
pressure is not smooth.
Table~\ref{tab:stiff-pressure-structure} shows that both
defects decrease as $\gamma$ increases, indicating an
approach to the complementarity relation and the interior
Hele--Shaw equation on the fixed grid.

\begin{table}[htbp]
  \centering
  \small
  \caption{Defects in the complementarity relation and
  the interior Hele--Shaw pressure equation
  for increasing $\gamma$.}
  \label{tab:stiff-pressure-structure}
  \vspace{3mm}
  \setlength{\tabcolsep}{3pt}
  \renewcommand{\arraystretch}{1.20}
  \begin{tabular}{c ccccccc}
    \hline
    $\gamma$
    & $10$ & $20$ & $40$ & $80$
    & $160$ & $320$ & $640$\\
    \hline
    $\mathcal C_\gamma$
    & $1.07\mathrm{E}{-1}$
    & $5.48\mathrm{E}{-2}$
    & $2.76\mathrm{E}{-2}$
    & $1.39\mathrm{E}{-2}$
    & $6.96\mathrm{E}{-3}$
    & $3.48\mathrm{E}{-3}$
    & $1.74\mathrm{E}{-3}$\\
    $\mathcal I_\gamma$
    & $1.91\mathrm{E}{-1}$
    & $9.24\mathrm{E}{-2}$
    & $4.52\mathrm{E}{-2}$
    & $2.26\mathrm{E}{-2}$
    & $1.09\mathrm{E}{-2}$
    & $4.96\mathrm{E}{-3}$
    & $2.84\mathrm{E}{-3}$\\
    \hline
  \end{tabular}
\end{table}

\FloatBarrier

\subsection{Two-dimensional tumor growth}
\label{subsec:num-2d-three-species}

Finally, we use the two-dimensional scheme in Appendix~\ref{app:2d}
to examine how the spatial distribution of species with
different growth rates affects tumor morphology.
The initial configuration is motivated by multispecies
tumor models with spatially heterogeneous compositions
\cite{FrieboesJinChuangWiseLowengrubCristini2010,WiseLowengrubFrieboesCristini2008}.

We take $\Omega=(-10,10)^2$, $N_x=N_y=500$,
$\tau=10^{-2}$ and $T=2.5$.
For three species without switching, the reaction matrix is
\begin{equation}
\label{eq:2d-reaction-matrix}
  \mathbf R(p)
  =
  \operatorname{diag}
  \bigl(8(1-p),\,4(1-p),\,1.5(1-p)\bigr).
\end{equation}

The initial total density is
\begin{equation}
\label{eq:2d-initial-total}
  n^0(x,y)
  =
  0.49
  \Big[
    1-
    \tanh\Big(
      \frac{\sqrt{x^2+y^2}-2}{0.08}
    \Big)
  \Big].
\end{equation}
To prescribe the initial composition, we set
\[
\begin{aligned}
  w_1(x,y)
  &=
  60\exp\Big(
    -\frac{(x-2.05)^2+y^2}{2(0.25)^2}
  \Big), \qquad  \quad
  w_2(x,y)
  =1,\\
  w_3(x,y)
  &=
  10\exp\Big(
    -\frac{(x+1)^2+(y+0.9)^2}{2(0.40)^2}
  \Big),
\end{aligned}
\]
and define
\begin{equation}
\label{eq:2d-initial-species}
  \phi_i^0(x,y)
  =
  \frac{w_i(x,y)}
       {\sum_{\ell=1}^3w_\ell(x,y)},
  \qquad
  n_i^0(x,y)
  =
  \phi_i^0(x,y)n^0(x,y),
  \qquad i=1,2,3.
\end{equation}
Initially, species 1 is concentrated near the right
edge of the tumor, species 2 occupies most of the
remaining region and species 3 lies mainly in the
lower-left interior.
Although the initial total density is radially symmetric,
the different growth rates and spatial distributions
of the species produce an asymmetric total growth term.

\begin{figure}[htbp]
  \centering
  \includegraphics[
    width=0.83\textwidth,
    trim=0 1.5cm 0 1.5cm,
    clip
  ]{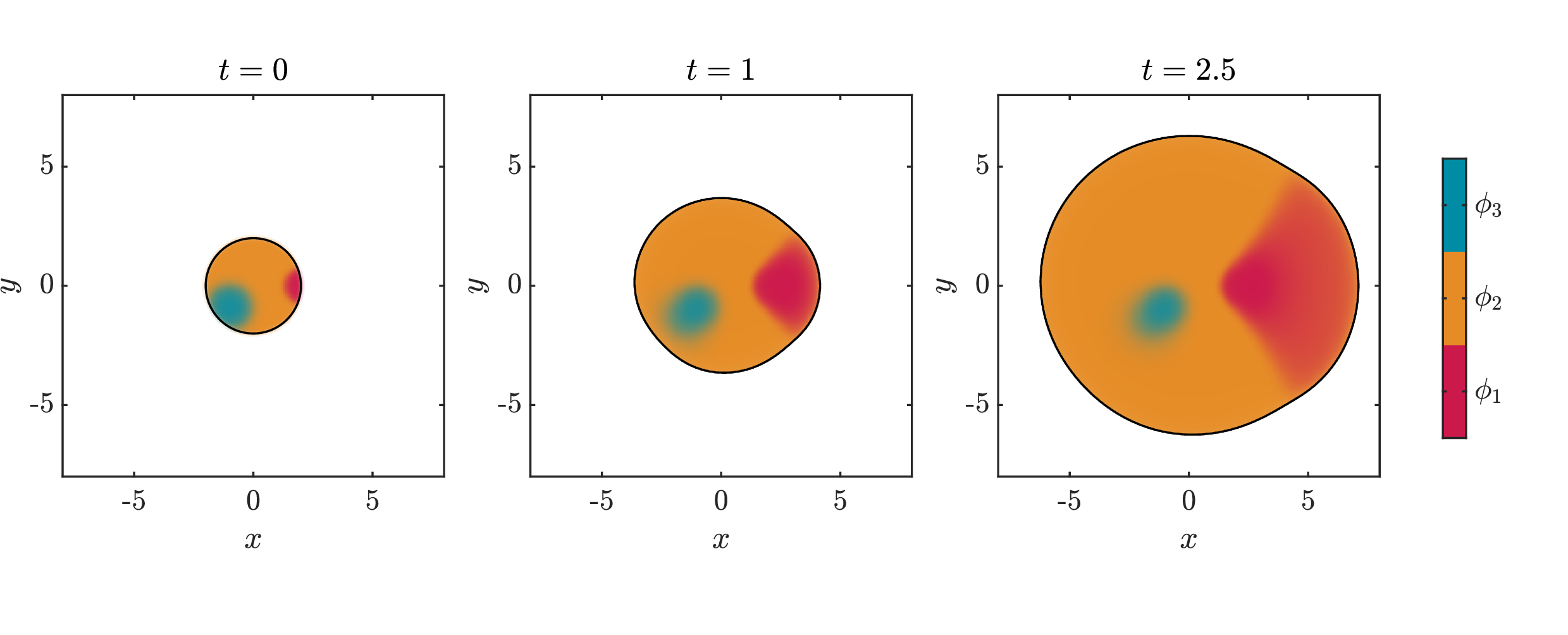}
  \includegraphics[
    width=0.83\textwidth,
    trim=0 2.3cm 0 1.5cm,
    clip
  ]{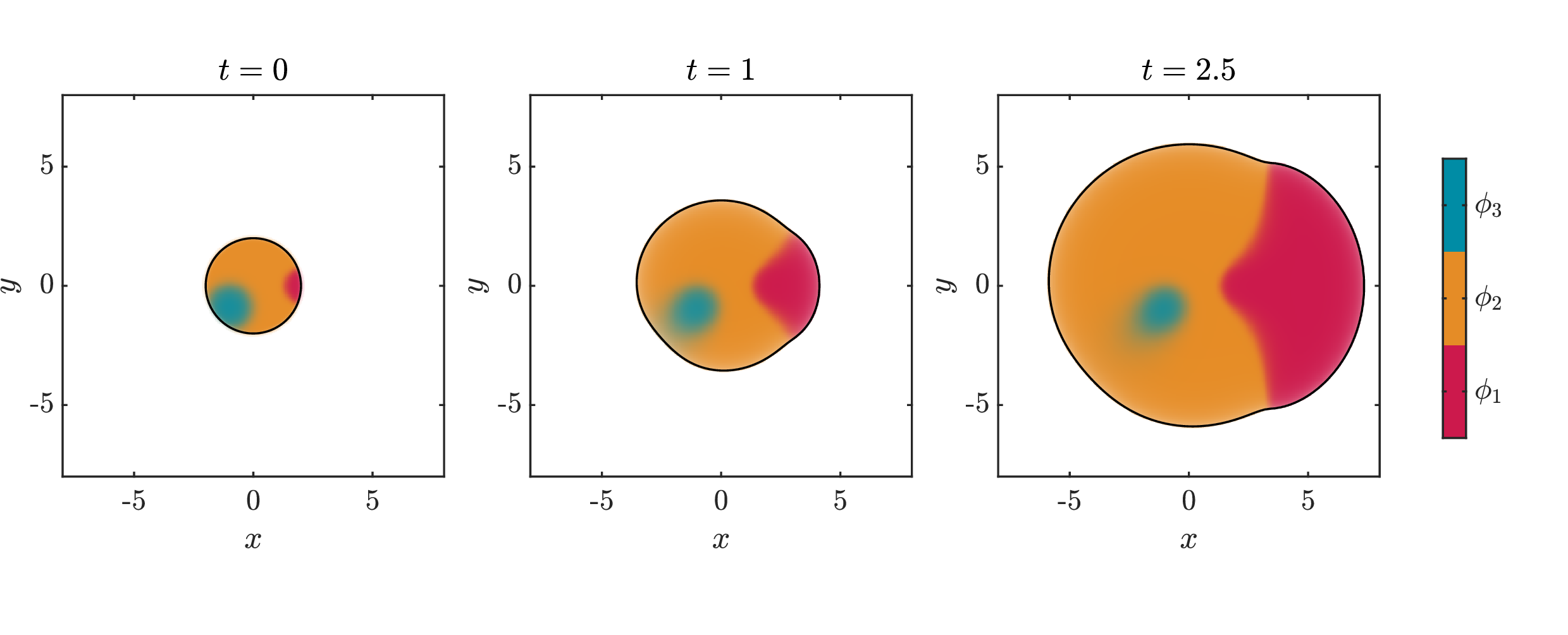}
\caption{
Evolution of the species composition
$\boldsymbol\phi=(\phi_1,\phi_2,\phi_3)^T$
for $\gamma=20$ (top) and $\gamma=3$ (bottom)
at $t=0,1,2.5$.
Colors indicate the species fractions, with opacity
reflecting the total density.
White denotes $n<0.02$, and black curves mark $n=0.5$.
}
  \label{fig:ex4-localized-clone}
\end{figure}

We compare the evolution for $\gamma=3$ and $\gamma=20$,
keeping the initial total density, species composition,
growth rates and numerical parameters unchanged.
Since $p^0=(n^0)^\gamma$, the initial pressure is lower
for $\gamma=20$ in the tumor interior, so growth is
initially less inhibited by pressure.
Changing $\gamma$ also affects the pressure gradient
and hence the Darcy transport.

Figure~\ref{fig:ex4-localized-clone} shows that, in both
cases, the fast-growing species spreads along the right
side of the moving front, while the slow-growing species
remains mainly in the interior.
At $T=2.5$, the $n=0.5$ contour is nearly circular
for $\gamma=20$, whereas it exhibits a more pronounced
rightward protrusion for $\gamma=3$.

\section{Conclusion}
\label{sec:conclusion}

We have developed an Onsager-based scheme for multispecies
pressure-driven tumor growth with a common homeostatic pressure.
The scalar Onsager formulation is coupled with compatible species
updates, so the global nonlinear problem remains scalar and the
species densities reconstruct the updated total density exactly.
The resulting first-order scheme is uniquely solvable and preserves
nonnegativity, the pressure bound and energy dissipation without
a time-step restriction.

The predictor--corrector extension improves temporal accuracy
while retaining these structural properties, with numerical tests
showing lower computational costs at comparable errors.
Extending this Onsager construction to species-dependent
homeostatic pressures is a direction for future work.

\bigskip
\noindent{\large\bf Acknowledgements.} X. Ruan acknowledges support  from the National Natural Science Foundation of China under grant  12201436 and the R\&D Program of Beijing Municipal Education Commission under grant KM202310028016.
W. Huang acknowledges support  from the National Natural Science Foundation of China under grant  12001034.

\appendix
\section{Two-dimensional extension}
\label{app:2d}

We extend the first-order scheme to a Cartesian mesh of
rectangular cells \(K_c\), with cell area \(|K_c|\).
For an interior edge \(e\) between cell $c$ and cell $d$, let \(|e|\) denote its
length and \(d_{cd}\) the distance between the centers
of \(K_c\) and \(K_d\).

\paragraph{Step 1: Scalar density step.}

The reaction mobility \(M_c^k\) is defined cellwise
as in the one-dimensional scheme.
For each interior edge \(e\) between cell $c$ and cell $d$, define
\[
  n_e^k=\frac{n_c^k+n_d^k}{2},
  \qquad
  J_{c,e}^{k+1}
  =
  -n_e^k
  \frac{p_d^{k+1}-p_c^{k+1}}{d_{cd}},
  \qquad
  J_{d,e}^{k+1}=-J_{c,e}^{k+1}.
\]
Here \(J_{c,e}^{k+1}\) is the outward normal flux
of the total density across \(e\), relative to \(K_c\).
We set \(J_{c,e}^{k+1}=0\) on boundary edges.

The scalar density step is
\begin{equation}
\label{eq:2d-density-step}
  n_c^{k+1}
  +
  \frac{\tau}{|K_c|}
  \sum_{e\subset\partial K_c}
  |e|J_{c,e}^{k+1}
  =
  \widehat n_c^{k+1},
  \qquad
  p_c^{k+1}=(n_c^{k+1})^\gamma,
\end{equation}
where
\[
  \widehat n_c^{k+1}
  =
  n_c^k-\tau M_c^k\bigl(p_c^{k+1}-p_H\bigr).
\]

\paragraph{Step 2: Species reaction step.}

We apply Step~2 of the one-dimensional scheme cellwise
to obtain the post-reaction species densities
\(\widehat{\boldsymbol n}_c^{k+1}\).
This update satisfies
\[
  \widehat{\boldsymbol n}_c^{k+1}\geq\boldsymbol0,
  \qquad
  \mathbf1^\top\widehat{\boldsymbol n}_c^{k+1}
  =
  \widehat n_c^{k+1}.
\]

\paragraph{Step 3: Species transport step.}

Let \(\boldsymbol\phi_c^{k+1}\) denote the unknown
composition vector in the transport step.
For each interior edge \(e=c|d\), define the upwind
species flux
\[
  \boldsymbol J_{c,e}^{k+1}
  =
  \bigl(J_{c,e}^{k+1}\bigr)^+
  \boldsymbol\phi_c^{k+1}
  +
  \bigl(J_{c,e}^{k+1}\bigr)^-
  \boldsymbol\phi_d^{k+1},
\]
and set \(\boldsymbol J_{c,e}^{k+1}=\boldsymbol0\)
on boundary edges.
The conservative transport update is
\begin{equation}
\label{eq:2d-species-conservative}
  n_c^{k+1}\boldsymbol\phi_c^{k+1}
  +
  \frac{\tau}{|K_c|}
  \sum_{e\subset\partial K_c}
  |e|\boldsymbol J_{c,e}^{k+1}
  =
  \widehat{\boldsymbol n}_c^{k+1}.
\end{equation}

Define
\[
  D_c^{k+1}
  =
  n_c^{k+1}
  +
  \frac{\tau}{|K_c|}
  \sum_{e\subset\partial K_c}
  |e|\bigl(J_{c,e}^{k+1}\bigr)^+.
\]
Then \eqref{eq:2d-species-conservative} can be written as
\[
  D_c^{k+1}\boldsymbol\phi_c^{k+1}
  =
  \widehat{\boldsymbol n}_c^{k+1}
  +
  \frac{\tau}{|K_c|}
  \sum_{\substack{e=c|d\\J_{c,e}^{k+1}<0}}
  |e|\bigl(-J_{c,e}^{k+1}\bigr)
  \boldsymbol\phi_d^{k+1}.
\]
If \(J_{c,e}^{k+1}<0\), then
\(p_d^{k+1}>p_c^{k+1}\), so every inflow comes
from a cell with strictly higher pressure and positive
density.
We therefore solve these equations in nonincreasing
pressure order.
Cells with equal pressure exchange no flux and may
be processed in any order.
For \(D_c^{k+1}>0\), we compute
\(\boldsymbol\phi_c^{k+1}\) from the equation above
and set
\(\boldsymbol n_c^{k+1}
=n_c^{k+1}\boldsymbol\phi_c^{k+1}\).
If \(D_c^{k+1}=0\), we set
\(\boldsymbol n_c^{k+1}=\boldsymbol0\).

The arguments for the one-dimensional scheme give
unique solvability, species nonnegativity, the pressure
bound, total-mass monotonicity and exact reconstruction
of the total density.
The compositions therefore belong to \(\Delta_m\)
on occupied cells.
The energy estimate in
Proposition~\ref{prop:multi-scalar-step} extends
with cell weights \(|K_c|\) and interior-edge weights
\(|e|/d_{cd}\), each interior edge being counted once.
These properties hold without a time-step restriction.

{\small
\bibliographystyle{siamplain-local}
\bibliography{references}
}

\end{document}